\documentclass[12pt]{amsart}
\usepackage{amssymb, amsmath, amscd, dcpic, pictexwd}

\newtheorem{thm}{Theorem}[section]
\newtheorem{qst}[thm]{Question}
\newtheorem{cor}[thm]{Corollary}
\newtheorem{lem}[thm]{Lemma}

\newtheorem{defn}[thm]{Definition}
\newtheorem{rem}[thm]{Remark}
\newtheorem{exm}[thm]{Example}
\numberwithin{equation}{section}
\newcommand{\ci}{\circ_{n}}

\newcommand{\be}{\beta}
\newcommand{\ga}{\gamma}
\newcommand{\Bx}{\beta^x}
\newcommand{\By}{\beta^y}
\newcommand{\Bh}{\beta^h}
\newcommand{\Gx}{\gamma^{x'}}
\newcommand{\Gy}{\gamma^{y'}}
\newcommand{\Gh}{\gamma^{h'}}

\newcommand{\tX}{\theta^x}
\newcommand{\tY}{\theta^y}
\newcommand{\tH}{\theta^h}

\newcommand{\vv}{\mathcal V}
\newcommand{\A}{\mathcal A}
\newcommand{\B}{\mathcal B}

\newcommand{\D}{\mathcal D}

\newcommand{\Se}{\mathcal S}

\newcommand{\LL}{\mathcal{L}}

\newcommand{\E}{\mathcal E}
\newcommand{\W}{\mathcal W}

\newcommand{\G}{\mathfrak{g}}
\newcommand{\U}{\mathfrak{U}}

\newcommand{\hh}{\mathfrak{h}}
\newcommand{\bb}{\mathfrak{b}}

\newcommand{\CC}{\mathbb{C}}
\newcommand{\ZZ}{\mathbb Z}

\newcommand{\One}{{\bf 1}}

\author{Bong H. Lian and Andrew R. Linshaw}
\address{Department of Mathematics\\ Brandeis University\\
Waltham, MA 02454.}

\title{Howe Pairs in the Theory of Vertex Algebras}
\begin{document}
\tableofcontents \maketitle
\begin{abstract}
For any vertex algebra $\vv$ and any subalgebra $\A\subset \vv$,
there is a new subalgebra of $\vv$ known as the {\it commutant} of
$\A$ in $\vv$. This construction was introduced by Frenkel-Zhu, and
is a generalization of an earlier construction due to Kac-Peterson
and Goddard-Kent-Olive known as the coset construction. In this
paper, we interpret the commutant as a vertex algebra notion of
invariant theory. We present an approach to describing commutant algebras in an appropriate category of vertex algebras by reducing the problem to a question in commutative algebra. We give an interesting example of a Howe pair (ie, a pair of mutual commutants) in the vertex algebra setting.

\end{abstract}

\section{Introduction}
For any vertex algebra $\vv$ and any subalgebra $\A\subset \vv$, the
commutant of $\A$ in $\vv$, denoted by $Com(\A,\vv)$, is defined to
be the set of vertex operators $v(z)\in \vv$ such that
$[a(z),v(w)]=0$ for all $a(z)\in\A$. This construction is analogous
to the ordinary commutant in the theory of associative algebras, and
was introduced by Frenkel-Zhu in \cite{FZ}, generalizing a previous
construction in representation theory \cite{KP} and conformal field
theory \cite{GKO} known as the coset construction. Describing
$Com(\A,\vv)$ by giving generators and OPE relations is generally a
non-trivial problem. A priori, it is far from clear when
$Com(\A,\vv)$ is finitely generated as a vertex algebra, even when
$\A$ and $\vv$ are finitely generated.

Equivalently, $Com(\A,\vv)$ is the subalgebra
$$\{v(z)\in\vv|\ a(z)\ci v(z)=0,\ \forall a(z)\in\A,\ n\ge 0\}.$$
Thus if we regard $\vv$ as a module over $\A$ via the \lq\lq left
regular action," $Com(\A,\vv)$ is the subalgebra of $\vv$ which is
annihilated by the operators $\{\hat{a}(n)|\ a\in\A,\ n\ge 0\}$. We
regard $\vv$ equipped with its $\A$-module structure as the analogue
of an associative algebra equipped with a Lie group or Lie algebra
action, and we regard $Com(\A,\vv)$, which we often denote by
$\vv^{\A_+}$, as the invariant subalgebra. Often $\A$ will be a homomorphic image of a current algebra $\mathcal{O}(\G,B)$, where $\G$ is a Lie algebra and
$B$ is a symmetric, invariant bilinear form on $\G$. In this case,
$\vv^{\A_+}$ is just the invariant space $\vv^{\G[t]}$, where
$\G[t]$ is the Lie subalgebra of the loop algebra
$\G[t,t^{-1}]$ generated by $\{ut^n|\ u\in \G,n\ge 0\}$. The problem of describing $\vv^{\G[t]}$ lies outside the realm of classical invariant theory since $\G[t]$ is both infinite-dimensional and
non-reductive. 

\subsection{Acknowledgements}

We thank G. Schwarz for many helpful discussions about classical invariant theory. After this work was completed, E. Frenkel communicated to us that in his study of the semi-infinite Weil complex of $sl(2)$, he had found (but not published) some of the same results in our paper. We thank him for bringing this to our attention.

\subsection{Howe pairs}
For any vertex algebra $\vv$ and subalgebra $\A\subset\vv$, we have $\A\subset Com(Com(\A,\vv),\vv))$. If this inclusion is an equality, so that $\A$ and $Com(\A,\vv)$ are mutual commutants, we say that $\A$ and $Com(\A,\vv)$ form a {\it Howe pair} inside $\vv$. Our main goal is to give an interesting example of a Howe pair in the vertex algebra setting, as well as outline a general strategy for describing commutant algebras of the form $\vv^{\G[t]}$ in an appropriate category of vertex algebras. 

We will focus on a particular situation which is induced by a problem in classical invariant theory. Let $\G$ be a finite-dimensional, semisimple, complex Lie algebra, and let $V$ be a finite-dimensional complex vector space which is a $\G$-module via
$\rho:\G\rightarrow End(V)$. Associated to $V$ is a vertex
algebra $\Se(V)$ known as a $\be\ga$-ghost system or a semi-infinite
symmetric algebra \cite{FMS}. The map $\rho$ induces a
vertex algebra homomorphism \begin{equation} \hat{\rho}:\mathcal{O}(\G,B)\rightarrow
\Se(V),\end{equation} where $B$ is the bilinear form $B(u,v)=
-Tr\big(\rho(u)\rho(v)\big)$. Letting $\Theta =
\hat{\rho}(\mathcal{O}\big(\G,B)\big)$, we will study the commutant
$\Se(V)^{\Theta_+}$. Generically, $\Se(V)^{\Theta_+}$ is a conformal
vertex algebra with conformal weight grading
$$\Se(V)^{\Theta_+}=\bigoplus_{n\ge 0}\Se(V)^{\Theta_+}_n,$$ and the weight-zero
subspace $\Se(V)^{\Theta_+}_0$ coincides with the classical ring
$Sym(V^*)^{\G}$ of invariant polynomial functions on $V$. In other
words, $\Se(V)^{\Theta_+}$ is a \lq\lq chiralization" of
$Sym(V^*)^{\G}$. 

\subsection{The Zhu functor and invariant differential operators}
In \cite{Z}, Zhu introduced a functorial construction which attaches
to every vertex algebra $\vv$ an associative algebra $A(\vv)$ known
as the {\it Zhu algebra} of $\vv$, together with a surjective linear
map $\pi_{Zh}:\vv\rightarrow A(\vv)$ known as the {\it Zhu map}. It
is well known that $A\big(\mathcal{O}(\G,B)\big)$ is the universal
enveloping algebra $\U\G$, and $A\big(\Se(V)\big)$ is the Weyl
algebra $ \D(V)$ of polynomial differential operators of $V$. $\D(V)$
has generators $\Bx,\Gx$ which are linear in $x\in V$, $x'\in V^*$,
and satisfy the commutation relations
\begin{equation}[\Bx,\Gx]=\langle x',x\rangle.\end{equation}
If we fix a basis $x_1,\dots,x_n$ for $V$ and a corresponding dual
basis $x'_1,\dots,x'_n$ for $V^*$, the variables $\ga^{x'_i}$
correspond to the linear functions $x'_i$, and the
variables $\be^{x_i}$ correspond to the first-order differential
operators $\frac{\partial}{\partial x'_i}$.

If $V$ is a $\G$-module via $\rho:\G\rightarrow End(V)$, there is an
induced action $\rho^*$ of $\G$ on $\D(V)$. We would like to study
the relationship between $\Se(V)^{\Theta_+}$ and the classical ring
$\D(V)^{\G}$ of invariant polynomial differential operators on $V$. Our discussion of $\D(V)^{\G}$ is based on \cite{S}.
The invariant subalgebra $\D(V)^{\G}$ contains $Sym(V^*)^{\G}$ as
the subspace of zeroth order invariant differential operators.
Recall that $\D(V)$ has a filtration known as the {\it Bernstein
filtration}
$$0\subset \D_0(V)\subset \D_1(V)\subset\cdots,$$ where $\sum
(\ga^{x'_i})^{n_i}(\be^{x_j})^{m_j}\in \D_n(V)$ iff $\sum_i n_i +
\sum_j m_j\le n$. It follows from (1.2) that the associated graded
object $$gr\big(\D(V)\big) = \bigoplus_{n>0}
\D_n(V)/\D_{n-1}(V)=Sym(V\oplus V^*).$$ Moreover, $\G$ acts on
$\D(V)$ by derivations of degree 0, so the above filtration
restricts to a filtration
$$0\subset \D_0(V)^{\G}\subset \D_1(V)^{\G}\subset\cdots$$ of
$\D(V)^{\G}$, and we have $$ gr\big(\D(V)^{\G}\big) =
gr\big(\D(V)\big)^{\G} = Sym(V\oplus V^*)^{\G},$$

The action of $\G$ on $\D(V)$ can be realized
by inner derivations. We have a Lie algebra homomorphism
$\tau:\G\rightarrow \D(V)$ given in our chosen basis by
\begin{equation}\tau(u) = -\sum_i
\be^{\rho(u)(x_i)}\ga^{x'_i},\end{equation} which we may extend to a
map $\U\G\rightarrow \D(V)$, and the action of $u\in\G$ on $\omega\in\D(V)$ is given by 
$\rho^*(u)(\omega) =
[\tau(u),\omega]$. Thus $\D(V)^{\G}$ may be
alternatively described as the commutant
$Com\big(T,\D(V)\big)$, where $T=\tau(\U\G)\subset \D(V)$,
in the sense of ordinary associative algebras. We have the following
commutative diagram:

\bigskip

\centerline{
\begindc{\commdiag}[3]
\obj(10,40)[SG]{$\Se(V)^{\Theta_+}$}
\obj(40,40)[S]{$\Se(V)$}
\obj(10,15)[DG]{$\D(V)^{\G}$}
\obj(40,15)[D]{$\D(V)$}
\mor{SG}{DG}{$\pi_{Zh}|_{\Se(V)^{\Theta_+}}$}
\mor{DG}{D}{$i$}[\atleft,\injectionarrow]
\mor{SG}{S}{$i$}[\atleft,\injectionarrow]
\mor{S}{D}{$\pi_{Zh}$}
\enddc}

\bigskip

The horizontal maps above are inclusions, and the vertical map on
the left is the restriction of the Zhu map on $\Se(V)$ to the
subspace $\Se(V)^{\Theta_+}$. A priori, this map
need not be surjective, and $\D(V)^{\G}$ need not coincide with the Zhu algebra of
$\Se(V)^{\Theta_+}$. Even when this map is surjective, so that any
set of generators $\omega_1,\dots,\omega_k$ of $\D(V)^{\G}$ lifts to
a set of vertex operators $\omega_1(z),\dots,\omega_k(z)$ in
$\Se(V)^{\Theta_+}$, it is not clear when this collection generates
$\Se(V)^{\Theta_+}$ as a vertex algebra.

For any $\G$-module $V$, $\D(V)^{\G}$ always contains the Euler
operator $\sum_i\be^{x_i}\ga^{x'_i}$, where $\{x_i\}$ is a basis of $V$
and $\{x'_i\}$ is the corresponding dual basis of $V^*$. If $V$ admits
a symmetric, $\G$-invariant bilinear form $B$, there is a Lie
algebra homomorphism $\psi: sl(2)\rightarrow \D(V)^{\G}$ given in an
orthonormal basis (relative to $B$) by the formulas
\begin{equation} h\mapsto  \sum_i\be^{x_i}\ga^{x'_i},\ \ \ \
x\mapsto \frac{1}{2} \sum_i \ga^{x'_i}\ga^{x'_i},\ \ \ \ y \mapsto
-\frac{1}{2} \sum_i \be^{x_i}\be^{x_i}.\end{equation} Here $x,y,h$ denote
the standard generators of $sl(2)$, satisfying $$[x,y]=h,\ \ \ \ [h,x]=2x,\ \ \ \ [h,y]=-2y.$$
$\psi$ may be extended to a map $\U(sl(2))\rightarrow \D(V)^{\G}$, and
we denote the image $\psi(\U(sl(2))\subset \D(V)^{\G}$ by  $A$.

Likewise in the vertex algebra setting, $\Se(V)^{\Theta_+}$ contains
a vertex operator analogous to the Euler operator above, which
generates a Heisenberg vertex algebra inside $\Se(V)^{\Theta_+}$ of
central charge $-dim V$. When $V$ is admits a symmetric, $\G$-invariant bilinear form, the map
$\psi:\U(sl(2))\rightarrow \D(V)^{\G}$ gives rise to a vertex
algebra homomorphism
\begin{equation}\hat{\psi}:\mathcal{O}\big(sl(2),-\frac{dim V}{8}K\big)\rightarrow
\Se(V)^{\Theta_+},\end{equation} where $K$ denotes the Killing form
on $sl(2)$. This construction is compatible with the Zhu
functor in the sense that the diagram below commutes:
\bigskip

\centerline{
\begindc{\commdiag}[3]
\obj(10,40)[OS]{$\mathcal{O}(sl(2),-\frac{dim V}{8}K)$}
\obj(40,40)[SG]{$\Se(V)^{\Theta_+}$}
\obj(10,15)[US]{$\U(sl(2))$}
\obj(40,15)[DG]{$\D(V)^{\G}$}
\mor{SG}{DG}{$\pi_{Zh}|_{\Se(V)^{\Theta_+}}$}
\mor{OS}{US}{$\pi_{Zh}$}
\mor{OS}{SG}{$\hat{\psi}$}
\mor{US}{DG}{$\psi$}
\enddc}

\bigskip

Let $\A$ denote the image of $\mathcal{O}\big(sl(2),-\frac{dim V}{8}K\big)$ under $\hat{\psi}$. Clearly $\pi_{Zh}(\A) = A$ since $\pi_{Zh}$ maps the generators of $\A$ to the generators of $A$. 

\subsection{Some open questions}
We regard $\Se(V)^{\Theta_+}$ as a vertex algebra analogue of the classical invariant ring $\D(V)^{\G}$, and we ask whether various properties of $\D(V)^{\G}$ have appropriate analogues in the vertex algebra setting. For example, $\D(V)^{\G}$ is finitely generated as a ring by a classical theorem of Hilbert \cite{E}. Working in $gr(\D(V)) = Sym(V\oplus V^*)$, the idea of the proof is to use the complete reducibility of the $\G$-action on $Sym(V\oplus V^*)$ to express $Sym(V\oplus V^*)^{\G}$ as a direct
summand. It is a standard fact in commutative
algebra that any ring which is a summand of a finitely generated polynomial ring is finitely generated \cite{E}. 

\begin{qst}Is $\Se(V)^{\Theta_+}$ finitely generated as a vertex algebra? Can one find a set of generators? Is this an appropriate analogue of Hilbert's theorem? More generally, when are commutant algebras of the form $\vv^{\G[t]}$ finitely generated?
\end{qst}

Unfortunately, $\Se(V)$ is not unitary as an $\mathcal{O}(\G,B)$-module in general, so a priori $\Se(V)$ need not decompose into a direct sum of irreducible $\mathcal{O}(\G,B)$-modules, and a similar proof cannot be expected to go through. One of our goals will be to outline an alternative approach to answering this kind of question.

Another classical question one can ask is whether $T=\tau(\U\G)$ and $\D(V)^{\G}$ form a Howe pair (i.e., a pair of mutual commutants) inside $\D(V)$. This question has been studied by Knop in \cite{Kn} in a much wider context, namely, when the linear space $V$ is replaced by an algebraic variety with an algebraic group action. 

\begin{qst} When do $\Theta$ and $\Se(V)^{\Theta_+}$ form a Howe pair inside $\Se(V)$?
\end{qst}

To answer Question 1.2, one needs to compute $Com(\Se(V)^{\Theta_+},\Se(V))$ and determine whether it coincides with $\Theta$. This may be possible to carry out even without a complete description of $\Se(V)^{\Theta_+}$. A priori, we have $\Theta\subset Com(\Se(V)^{\Theta_+},\Se(V))$. Note that if $\B$ is any subalgebra of $\Se(V)^{\Theta_+}$, we have
$$Com(\Se(V)^{\Theta_+},\Se(V))\subset Com(\B,\Se(V)).$$ If we can show that $Com(\B,\Se(V)) = \Theta$, it follows that
$$\Theta\subset Com(\Se(V)^{\Theta_+},\Se(V))\subset Com(\B,\Se(V)) = \Theta,$$ so all of these algebras are equal.
\subsection{Statement of main result}
We will answer Question 1.2 in the basic but non-trivial special case when $\G = sl(2)$ and $V$ is the adjoint module. We will show that $\Se(V)^{\A_+} = \Theta$, where $\A$ is the subalgebra $\hat{\psi}\big(\mathcal{O}\big(sl(2),-\frac{3}{8}K\big)\big)$ of $\Se(V)^{\Theta_+}$, as above. It follows immediately that $Com(\Se(V)^{\Theta_+},\Se(V)) = \Theta$. Thus we obtain


\begin{thm} In the case $\G=sl(2)=V$, the subalgebras $\Theta$ and $\Se(V)^{\Theta_+}$ form a Howe pair inside $\Se(V)$.
\end{thm}
Note that both $\Se(V)^{\A_+}$ and $\Se(V)^{\Theta_+}$ are commutant algebras of the form $\vv^{\G[t]}$. In the case $\G=sl(2)=V$, $\Se(V)^{\A_+}$ is indeed finitely generated; it is just a copy of $\mathcal{O}(sl(2),-K)$. We expect that our method for calculating $\Se(V)^{\A_+}$ in this special case will useful for describing $\Se(V)^{\Theta_+}$ and $\Se(V)^{\A_+}$ for more general $\G$ and $V$, and possibly more general  commutant algebras of the form $\vv^{\G[t]}$ as well. We hope to return to these questions in the future. 
 
\subsection{Outline of proof} 
Following ideas introduced in \cite{L}, we reduce the problem of computing $\Se(V)^{\A_+}$ to a question in commutative algebra. We will single out a
certain category $\Re$ of $\mathbb{Z}_{\ge 0}$-filtered vertex
algebras whose associated graded objects are $\mathbb{Z}_{\ge
0}$-graded supercommutative rings. In particular, the assignment
$\vv\mapsto gr(\vv)$ is a functor from $\Re$ to the category of
$\mathbb{Z}_{\ge 0}$-graded supercommutative rings. The filtrations possessed by the objects of $\Re$ are examples of the {\it good increasing filtrations} introduced by Li in \cite{Li2}. $\Re$ includes
all vertex algebras of the form $\Se(V)$ and $\mathcal{O}(\G,B)$,
and is closed under taking subalgebras, so $\Theta$, $\A$, and $\Se(V)^{\Theta_+}$ and $\Se(V)^{\A_+}$ lie in $\Re$ as well. Moreover, $\Re$ has an
important {\it reconstruction property}; if we can find a set of
generators for the ring $gr(\vv)$, we can use them to construct a
set of generators for $\vv$ as a vertex algebra.

In the case $\vv=\Se(V)$, the ring $gr\big(\Se(V)\big)$ is
isomorphic to the polynomial algebra
$$P=Sym\big(\bigoplus_{k\ge 0} (V_k\oplus V^*_k)\big),$$ where
each $V_k$ and $V^*_k$ are copies of $V$ and $V^*$, respectively.
The action of $\A = \hat{\psi}\big(\mathcal{O}\big(sl(2),-\frac{3}{8}K\big)\big)$ on $\Se(V)$ induces an action of
the Lie algebra $sl(2)[t]$ on $P$ by derivations of degree 0, and we
denote the $sl(2)[t]$-invariant subalgebra by $P^{\A_+}$. We will
study $\Se(V)^{\A_+}$ indirectly by studying its associated
graded algebra $gr\big(\Se(V)^{\A_+}\big)$, and comparing it to
$P^{\A_+}$. There is a canonical injective ring homomorphism
\begin{equation} \Gamma: gr\big(\Se(V)^{\A_+}\big)\hookrightarrow
P^{\A_+}.\end{equation} 

Using tools from commutative algebra and classical invariant theory, we will be able to write down generators for $P^{\A_+}$ and see that (1.6) is in fact an isomorphism. By the reconstruction property, we obtain generators for $\Se(V)^{\A_+}$ as a vertex algebra, and we will see explicitly that these generators coincide with the generators of $\Theta$.

\subsection{Related questions} $\Se(V)^{\Theta_+}$ is an interesting vertex algebra that appears in several other contexts as well. Associated to the vector space $V$ is another
vertex algebra $\E(V)$ known as a $bc$-ghost system or a
semi-infinite exterior algebra, which is an odd
analogue of $\Se(V)$ \cite{FMS}. If $V$ is a $\G$-module via
$\rho:\G\rightarrow End(V)$, there is an induced vertex algebra map
$\mathcal{O}(\G,B)\rightarrow \E(V)$, analogous to (1.1), where $B(u,v) =
Tr\big(\rho(u)\rho(v)\big)$. In the case $V = \G$, the tensor
product $\E(\G)\otimes \Se(\G)$ is known as the {\it semi-infinite
Weil complex} of $\G$ \cite{FF}. $\W(\G)$ is a conformal vertex algebra with weight grading by
the non-negative integers, and $\W(\G)_0$ coincides with the
classical Weil algebra $W(\G) = \Lambda(\G^*)\otimes Sym(\G^*)$.
$\W(\G)$ also has a $\mathbb{Z}$-grading by {\it fermion number},
and contains a BRST current $J(z)$ whose zero mode $J(0)$ is a
square-zero derivation of degree 1 in this grading. The complex
$\big(\W(\G)^*,J(0)\big)$ coincides with a certain {\it relative semi-infinite
cohomology complex} of the affine Lie algebra $\hat{\G}$ of central
charge $-1$, with coefficients in the module $\Se(\G)$ \cite{F}\cite{FGZ}. The
cohomology $H^*\big(\W(\G),J(0)\big)$ is analogous to the Lie
algebra cohomology of $\G$ with coefficients in $Sym(\G^*)$, and was
studied in \cite{FF} and \cite{AK}. It is related to the commutant
$\Se(\G)^{\Theta_+}$ since $\Se(\G)^{\Theta_+} =
Ker\big(J(0)\big)\bigcap \big(1\otimes\Se(\G)\big)$. In the case
$\G=sl(2)$, Akman wrote down several examples of vertex
operators in $\W(\G)$ which represent nonzero cohomology classes,
and her list includes the generators of the subalgebra $\A\subset \Se(\G)^{\Theta_+}$, which plays an important role in this paper.

The semi-infinite Weil complex can also be used to define a vertex
algebra valued equivariant cohomology theory for any smooth
$G$-manifold $M$, where $G$ is a compact Lie group, known as the
{\it chiral equivariant cohomology} \cite{LL}. The definition of
${\bf H}^*_G(M)$ is analogous to the de Rham model for the classical
equivariant cohomology $H^*_G(M)$ due to H. Cartan \cite{C1}\cite{C2}
and further developed in \cite{DKV}\cite{GS}. ${\bf H}^*_G(M)$ is
$\mathbb{Z}_{\ge 0}$-graded by weight, and contains $H^*_G(M)$ as the
weight-zero subspace. Taking $\G$ to be the complexified Lie algebra
of $G$, $\W(\G)$ plays the role of $W(\G)$ in the classical theory,
and the commutant construction plays the role of classical invariant
theory in defining the appropriate notion of basic subcomplex $\W(\G)_{bas}$.

When $G$ is simple and $M$ is a point, ${\bf H}^*_G(pt) = H^*(\W(\G)_{bas})$ is an
interesting conformal vertex algebra containing $H^*_G(pt) =
Sym(\G^*)^{G}$ as the weight-zero subspace. Computing ${\bf
H}^*_G(pt)$ is a fundamental building block of this theory since for
any $G$-manifold $M$, ${\bf H}^*_G(M)$ is a module over ${\bf
H}^*_G(pt)$ via a chiral analogue of the Chern-Weil map. Moreover, $\Se(\G)^{\Theta_+}$ is a canonical subalgebra of $\W(\G)_{bas}$, and we expect that describing $\Se(\G)^{\Theta_+}$ will be a key
step in computing ${\bf H}^*_G(pt)$.

\section{Vertex algebras}
In this section, we define vertex algebras and their modules, which
have been discussed from various different points of view in the
literature
\cite{B}\cite{FHL}\cite{FLM}\cite{FBZ}\cite{FMS}\cite{K}\cite{Li}\cite{LZ}\cite{MS}.
We will follow the formalism developed in \cite{LZ} and partly in
\cite{Li}. Let $V=V_0\oplus V_1$ be a super vector space over
$\mathbb{C}$, and let $z,w$ be formal variables. By $QO(V)$, we mean
the space of all linear maps
$$V\rightarrow V((z))=\{\sum_{n\in\mathbb{Z}} v(n) z^{-n-1}|
v(n)\in V,\ v(n)=0\ for\ n>>0 \}.$$ Each element $a\in QO(V)$ can be
uniquely represented as a power series
$a(z)=\sum_{n\in\mathbb{Z}}a(n)z^{-n-1}\in (End\ V)[[z,z^{-1}]]$,
although the latter space is clearly much larger than $QO(V)$. We
refer to $a(n)$ as the $n$-th Fourier mode of $a(z)$. Each $a\in
QO(V)$ is assumed to be of the shape $a=a_0+a_1$ where
$a_i:V_j\rightarrow V_{i+j}((z))$ for $i,j\in\ZZ/2$, and we write
$|a_i| = i$.

On $QO(V)$ there is a set of non-associative bilinear operations,
$\circ_n$, indexed by $n\in\ZZ$, which we call the $n$-th circle
products. For homogeneous $a,b\in QO(V)$ they are defined by
$$a(w)\circ_n b(w)=Res_z a(z)b(w)~\iota_{|z|>|w|}(z-w)^n-
(-1)^{|a||b|}Res_z b(w)a(z)~\iota_{|w|>|z|}(z-w)^n.$$
Here $\iota_{|z|>|w|}f(z,w)\in\mathbb{C}[[z,z^{-1},w,w^{-1}]]$ denotes the
power series expansion of a rational function $f$ in the region
$|z|>|w|$. Note that $\iota_{|z|>|w|}(z-w)^n \neq \iota_{|w|>|z|}(z-w)^n$ for $n<0$. We usually omit the symbol $\iota_{|z|>|w|}$ and just
write $(z-w)^{n}$ to mean the expansion in the region $|z|>|w|$,
and write $(-1)^n(w-z)^{n}$ to mean the expansion in $|w|>|z|$. It is
easy to check that $a(w)\circ_n b(w)$ above is a well-defined
element of $QO(V)$.

The non-negative circle products are connected through the {\it
operator product expansion} (OPE) formula (\cite{LZ}, Prop. 2.3).
For homogeneous $a,b\in QO(V)$, we have \begin{equation} a(z)b(w)=\sum_{n\ge 0}a(w)\circ_n
b(w)~(z-w)^{-n-1}+:a(z)b(w): \end{equation} as formal power series in $z,w$. Here
$$:a(z)b(w):\ =a(z)_-b(w)\ +\ (-1)^{|a||b|} b(w)a(z)_+\ ,$$ where
$a(z)_-=\sum_{n<0}a(n)z^{-n-1}$ and $a(z)_+=\sum_{n\ge
0}a(n)z^{-n-1}$. (2.1) is customarily written as
$$a(z)b(w)\sim\sum_{n\ge 0}a(w)\circ_n b(w)~(z-w)^{-n-1},$$ where
$\sim$ means equal modulo the term $:a(z)b(w):\ $.

Note that $:a(z)b(z):$ is a well-defined element of
$QO(V)$. It is called the {\it Wick product} of $a$ and $b$, and it
coincides with $a(z)\circ_{-1}b(z)$. The other negative circle products
are related to this by
$$ n!~a(z)\circ_{-n-1}b(z)=\ :(\partial^n a(z))b(z):\ ,$$
where $\partial$ denotes the formal differentiation operator
$\frac{d}{dz}$. For $a_1(z),...,a_k(z)\in QO(V)$, the $k$-fold
iterated Wick product is defined to be
$$:a_1(z)a_2(z)\cdots a_k(z):\ =\ :a_1(z)b(z):$$
where $b(z)=\ :a_2(z)\cdots a_k(z):\ $.

The set $QO(V)$ is a nonassociative algebra with the operations
$\circ_n$ and a unit $1$. We have $1\circ_n a=\delta_{n,-1}a$ for
all $n$, and $a\circ_n 1=\delta_{n,-1}a$ for $n\ge -1$. We are
interested in subalgebras $\A\subset QO(V)$, that is, linear subspaces
of $QO(V)$ containing 1, which are closed under the circle products.
In particular $\A$ is closed under formal differentiation $\partial$
since $\partial a=a\circ_{-2}1$. We call such a subalgebra a {\it circle algebra}
(also called a quantum operator algebra in \cite{LZ}). Many formal algebraic
notions are immediately clear: a circle algebra homomorphism is just a linear
map which sends $1$ to $1$ and preserves all circle products; a module over $\A$ is a
vector space $M$ equipped with a circle algebra homomorphism $\A\rightarrow
QO(M)$, etc. A subset $S=\{a_i|\ i\in I\}$ of $\A$ is
said to {\it generate} $\A$ if any element $a\in\A$ can be written
as a linear combination of nonassociative words in the letters
$a_i$, $\ci$, for $i\in I$ and $n\in\mathbb Z$. We say that $S$ {\it
strongly generates} $\A$ if any $a\in\A$ can be written as linear
combination of words in the letters $a_i,\ci$, for $n<0$.
Equivalently, $\A$ is spanned by the collection of vertex operators
of the form
$:\partial^{k_1}a_{i_1}(z)\cdots\partial^{k_m}a_{i_m}(z):$, for
$k_1,\dots,k_m \ge 0$.

\begin{rem}
Fix a nonzero vector $\One\in V$ and let $a,b\in QO(V)$ such that
$a(z)_+\One=b(z)_+\One=0$ for $n\ge 0$. Then it follows immediately
from the definition of the circle products that $(a\circ_p
b)_+(z)\One=0$ for all $p$. Thus if a circle algebra $\A$ is generated by
elements $a(z)$ with the property that $a(z)_+\One=0$, then every
element in $\A$ has this property. In this case the vector $\One$
determines a linear map
$$
\chi:\A\rightarrow V,~~~a\mapsto a(-1)\One=\lim_{z\rightarrow
0}a(z)\One$$ (called the creation map in \cite{LZ}), having the
following basic properties: \begin{equation}\chi(1)=\One,~~~\chi(a\circ_n b)=a(n)b(-1)\One,~~~\chi(\partial^p
a)=p! ~a(-p-1)\One.\end{equation}
\end{rem}

Next, we define the notion of commutativity in a circle algebra.
\begin{defn}We say that $a,b\in QO(V)$ circle commute if \begin{equation} (z-w)^N
[a(z),b(w)]=0\end{equation} for some $N\ge 0$. Here $[,]$ denotes the supercommutator. If $N$ can be chosen to be 0,
then we say that $a,b$ commute. A circle algebra is said to be
commutative if its elements pairwise circle commute.\end{defn}

Note that this condition implies that $a\ci b = 0$ for $n \ge N$. An
easy calculation gives the following very useful characterization of
circle commutativity.

\begin{lem} The condition (2.3) is equivalent to
the condition that the following two equations hold:
\begin{equation}[a(z)_+,b(w)]= \sum_{p=0}^{N-1} (a\circ_p
b)(w)(z-w)^{-p-1},\end{equation}
\begin{equation}[a(z)_-,b(w)]= \sum_{p=0}^{N-1} (-1)^p(a\circ_p
b)(w)(w-z)^{-p-1}.\end{equation}\end{lem}

The notion of a commutative circle algebra is abstractly equivalent
to the notion of a vertex algebra (see for example \cite{FLM}). Briefly, every commutative
circle algebra $\A$ is itself a faithful $\A$-module, called the
{\it left regular module}. Define
$$\rho:\A\rightarrow QO(\A),\ \ \ \ a\mapsto\hat a,\ \ \ \ \hat
a(\zeta)b=\sum_{n\in\mathbb{Z}} (a\circ_n b)~\zeta^{-n-1}.$$ It
can be shown (see \cite{L1} and \cite{LL}) that $\rho$ is an injective circle algebra
homomorphism, and the quadruple of structures
$(\A,\rho,1,\partial)$ is a vertex algebra in the sense of
\cite{FLM}. Conversely, if $(V,Y,\One,D)$ is a vertex algebra, the
collection $Y(V)\subset QO(V)$ is a commutative circle algebra.
{\it We will refer to a commutative circle algebra simply as a
vertex algebra throughout the rest of this paper}.

\begin{rem}
Let $\A'$ be the vertex algebra generated by $\rho(\A)$ inside
$QO(\A)$. Since $\hat a(n)1 = a(z)\ci 1 = 0$ for all $a\in\A$ and
$n\ge 0$, it follows from Remark 2.1 that for every $\alpha\in \A'$,
we have $\alpha_+ 1 = 0$. The creation map $\chi: \A'\rightarrow \A$
sending $\alpha\mapsto \alpha(-1)1$ is clearly a linear isomorphism
since $\chi\circ\rho = id$. It is often convenient to pass between
$\A$ and its image $\A'$ in $QO(\A)$. For example, we shall often
denote the Fourier mode $\hat{a}(n)$ simply by $a(n)$. When we say
that a vertex operator $b(z)$ is annihilated by the Fourier mode
$a(n)$ of a vertex operator $a(z)$, we mean that $a\ci b = 0$. Here
we are regarding $b$ as an element of the state space $\A$, while
$a$ operates on the state space, and the map $a\mapsto \hat{a}$ is
the state-operator correspondence.
\end{rem}

Let $\A$ be a vertex algebra, and let $a(z),b(z),c(z)\in \A$. The
following well-known formulas will be useful to us.

\begin{equation}:(:ab:)c:\  = \ :abc:
+ \sum_{n\ge 0}\frac{1}{(n+1)!}\big( :(\partial^{n+1} a)(b\ci c):\ +
(-1)^{|a||b|} (\partial^{n+1} b)(a\ci c):\big)\ .\end{equation}

For any $n\ge 0$, we have
\begin{equation} a\ci
(:bc:) -\ :(a\ci b)c:\ - (-1)^{|a||b|}\ :b(a\ci c): = \sum_{i=1}^n
\binom ni (a\circ_{n-i}b)\circ_{i-1}c.
\end{equation}
For any $n\in\mathbb{Z}$, we have
\begin{equation} a\ci b = \sum_{p\in\mathbb{Z}}
(-1)^{p+1}\big(b\circ_p a\big)\circ_{n-p-1}1.
\end{equation}

By the preceding remark, in order to prove these identities, it suffices to show that
$\hat{a},\hat{b},\hat{c}$ satisfy them, which can be
checked by applying the creation map to both sides and then using
(2.2). Equations (2.6), (2.7), and (2.8) measure the non-associativity of the Wick product, the failure of the positive circle products to be derivations of the Wick product, and the failure of
the circle products to be commutative, respectively.

\subsection{Conformal structure} Many vertex algebras $\vv$ have a
{\it Virasoro element}, that is, a vertex operator $L(z)$ satisfying
the OPE
\begin{equation} L(z)L(w) \sim
\frac{c}{2}(z-w)^{-4} + 2L(w)(z-w)^{-2} + \partial L(w)(z-w)^{-1},
\end{equation} where the constant $c$ is called the {\it central charge} of $L(z)$. It is customary to write $L(z) = \sum_{n\in\mathbb{Z}} L(n) z^{-n-1}$ in the form $\sum_{n\in\mathbb{Z}} L_n z^{-n-2}$, so that $L(n) = L_{n-1}$. Often
we require further that $L_0$ be diagonalizable and $L_{-1}$ acts on
$\vv$ by formal differentiation. In this case, the pair
$\big(\vv,L(z)\big)$ is called a {\it conformal vertex algebra of
central charge $c$}. An element $a(z)\in \vv$ which is an
eigenvector of $L_0$ with eigenvalue $\Delta$ is said to have {\it
conformal weight} $\Delta$, and we denote the subspace of conformal
weight $\Delta$ by $\vv_{\Delta}$. If $a(z)\in\vv_{\Delta}$
satisfies the OPE
$$L(z)a(w) \sim \Delta a(w)(z-w)^{-2} + \partial
a(w)(z-w)^{-1},$$ so that all the higher poles vanish, $a(z)$ is
said to be {\it primary}. In any conformal vertex algebra $\vv$, the
operator $\ci$ is homogeneous of weight $-n-1$. In particular, the
Wick product $\circ_{-1}$ is homogeneous of weight zero, so $\vv_0$
is closed under the Wick product. If the conformal weight grading is
a $\mathbb{Z}_{\ge 0}$-grading $\vv = \bigoplus_{n\ge 0} \vv_n$
(which will be the case in all our examples), $\vv_0$ is an
associative, supercommutative algebra with unit 1 under the Wick
product.

\begin{exm} {\it Current Algebras}
\end{exm}
Let $\G$ be a Lie algebra equipped with a symmetric $\G$-invariant
bilinear form $B$. The loop algebra of $\G$ is defined to be
$$\G[t,t^{-1}] = \G\otimes \mathbb C[t,t^{-1}],$$
with bracket given by $$[u t^n,v t^m]=[u,v] t^{n+m}.$$ The form $B$
determines a 1-dimensional central extension of $\G[t,t^{-1}]$
$$\hat{\G}= \G[t,t^{-1}]\oplus \mathbb{C}\kappa,$$
with bracket $$ [u t^n,vt^m]=[u,v] t^{n+m} + n B(u,v)
\delta_{n+m,0}\kappa.$$ $\hat{\G}$ is equipped with the
$\mathbb{Z}$-grading $deg(ut^n)=n$, and $deg(\kappa)=0$. Let
$\hat{\G}_{\ge 0} \subset \hat{\G}$ be the subalgebra of elements of
non-negative degree, and let $$N(\G,B) = \U\hat{\G}\otimes_{\hat{\G}_{\ge
0}}\bf{C}$$ be the induced $\hat{\G}$-module, where $\bf{C}$ is the
1-dimensional $\hat{\G}_{\ge 0}$-module on which $\G[t]$ acts by zero and
$\kappa$ acts by 1. Clearly $N(\G,B)$ is graded by the non-positive
integers. For each $u\in\G$, let $u(n)$ denote the linear operator
on $N(\G,B)$ representing $ut^n$, and put

\begin{equation}u(z)=\sum_{n\in \mathbb{Z}} u(n)z^{-n-1}\in QO\big(N(\G,B)\big).\end{equation}

The collection $\{u(z)|\ u\in\G\}$ generates a vertex algebra inside
$QO\big(N(\G,B)\big)$, which we denote by $\mathcal{O}(\G,B)$
\cite{FZ}\cite{LZ}\cite{L2}. For any $u,v\in\G$, the vertex
operators $u(z),v(z)\in \mathcal{O}(\G,B)$ satisfy the OPE
\begin{equation}u(z)v(w)\sim B(u,v)(z-w)^{-2}+[u,v](w)(z-w)^{-1}.
\end{equation}

Let ${\bf 1}$ denote the vacuum vector $1\otimes 1\in N(\G,B)$.
\begin{lem} \cite{L2} The creation map $\chi: \mathcal{O}(\G,B)\rightarrow
N(\G,B)$ sending $a(z)\rightarrow a(-1){\bf 1}$ is an isomorphism of
$\mathcal{O}(\G,B)$-modules.
\end{lem}
In fact, for $u_1,\dots,u_k\in\G$ and $n_1,\dots,n_k\ge 0$,
$$\chi\big(:\partial^{n_1}u_1(z)\cdots \partial^{n_k}u_k(z):\big)
= n_1! \cdots n_k! u_1(-n_1-1)\cdots u_k(-n_k-1).$$ By the
Poincare-Birkhoff-Witt (PBW) theorem, we may choose a basis of
$N(\G,B)$ consisting of monomials of the form $u_1(-n_1-1)\cdots
u_k(-n_k-1)$. Hence $\mathcal{O}(\G,B)$ is spanned by the collection
of standard monomials
\begin{equation} :\partial^{n_1}u_1(z)\cdots
\partial^{n_k}u_k(z):\ .
\end{equation}

If $\G$ is finite-dimensional and the form $B$ is non-degenerate,
$\mathcal{O}(\G,\lambda B)$ admits a Virasoro element $L_{\mathcal{O}}(z)$
such that $\big(\mathcal{O}(\G,\lambda B),L_{\mathcal{O}}(z)\big)$ is a
conformal vertex algebra, for all but finitely many values of
$\lambda\in\mathbb{C}$ \cite{L2}. For example, if $\G$ is simple,
$\mathcal{O}(\G,\lambda K)$ has a Virasoro element given by the {\it
Sugawara-Sommerfield formula}:
\begin{equation}L_{\mathcal{O}}(z) =
\frac{1}{2\lambda+1}\sum_i :u_i(z)u_i(z):\ ,\end{equation} whenever $\lambda \neq
-\frac{1}{2}$, where the $u_i$ form an orthonormal basis of $\G$
relative to the Killing form $K$. Note that we have chosen a
normalization so that we do not need to mention the dual Coxeter
number of $\G$. $L_{\mathcal{O}}(z)$ has central charge $\frac{2\lambda
dim(\G)}{2\lambda+1}$, and for each $u\in\G$, $u(z)$ is primary of
conformal weight 1. In fact, $L_{\mathcal{O}}(z)$ is characterized
by these properties \cite{L2}.

\begin{exm} {\it $\beta\gamma$-ghost systems}
\end{exm}
Let $V$ be a finite-dimensional vector space. Regard $V\oplus V^*$
as an abelian Lie algebra. Then its loop algebra has a
one-dimensional central extension
$$\mathfrak{h} = \hh(V) = (V\oplus V^*)[t,t^{-1}]\oplus \mathbb{C}\tau,$$
which is known as a Heisenberg algebra. Its bracket is given by
$$[(x,x')t^n,(y,y')t^m]=(\langle y',x\rangle-\langle
x',y\rangle)\delta_{n+m,0}\tau,$$ for $x,y\in V$ and $x',y'\in V^*$. Let $\bb\subset\hh$ be the subalgebra generated by $\tau$, $(x,0)t^n$, and $(0,x')t^m$, for $n\ge 0$ and $m>0$, and let $\bf{C}$
be the one-dimensional $\bb$-module on which $(x,0)t^n$ and
$(0,x')t^m$ act trivially and the central element $\tau$ acts by
the identity. Denote the linear operators representing
$(x,0)t^n,(0,x')t^n$ on $\U\hh\otimes_{\U\bb}\bf{C}$ by
$\beta^x(n),\gamma^{x'}(n-1)$, respectively, for $n\in\mathbb{Z}$.
The power series
$$\beta^x(z)=\sum_{n\in\mathbb{Z}}\beta^x(n)z^{-n-1},\ \ \ 
\gamma^{x'}(z)=\sum_{n\in\mathbb{Z}}\gamma^{x'}(n)z^{-n-1}\ \in
QO(\U\hh\otimes_{\U\bb}\bf{C})$$ generate a vertex algebra $\Se(V)$
inside $QO(\U\hh\otimes_{\U\bb}\bf{C})$, and the generators satisfy the OPE relations
\begin{equation}\beta^x(z)\gamma^{x'}(w)\sim\langle x',x\rangle (z-w)^{-1},\ \ \ \beta^x(z)\beta^y(w)\sim 0,\ \ \ \gamma^{x'}(z)\gamma^{y'}(w)\sim 0.\end{equation} This algebra was introduced in
\cite{FMS}, and is known as a $\be\ga$-ghost system, or a
semi-infinite symmetric algebra. By the PBW theorem, the
vector space $\U\hh\otimes_{\U\bb}\bf{C}$ has the structure of a
polynomial algebra with generators given by the negative Fourier
modes $\beta^x(n),\gamma^{x'}(n)$, $n<0$, which are linear in $x\in
V$ and $x'\in V^*$. It follows that $\Se(V)$ is spanned by the
collection of iterated Wick products of the form
$$\mu=\ :\partial^{n_1}\be^{x_1}\cdots
\partial^{n_s}\be^{x_s}
\partial^{m_1}\ga^{x'_1}\cdots \partial^{m_t}\ga^{x'_t}:\ .$$

$\Se(V)$ has the following Virasoro element:
\begin{equation}L_{\Se}(z) = \sum_i :\be^{x_i}(z)\partial\ga^{x'_i}(z):\ ,\end{equation}
where $x_1,\dots,x_n$ is a basis of $V$ and $x'_1,\dots,x'_n$ is the
corresponding dual basis of $V^*$. $L_{\Se}(z)$ is characterized by the property
that it is a Virasoro element of central charge $2 dim(V)$, and
$\Bx(z),\Gx(z)$ are primary of conformal weights $1,0$, respectively. 

Suppose that $V$ is a finite-dimensional $\G$-module via
$\rho:\G\rightarrow End(V)$, where $\G$ is a finite-dimensional Lie
algebra.

\begin{lem} The map $\rho$ induces a vertex algebra homomorphism
$$\hat{\rho}:\mathcal{O}(\G,B)\rightarrow\Se(V),$$ where $B$ is the
bilinear form $B(u,v) = -Tr\big(\rho(u)\rho(v)\big)$ on $\G$.
\end{lem}
\begin{proof} In terms of a basis $x_1,\dots,x_n$ for $V$ and dual basis
$x'_1,\dots,x'_n$ for $V^*$, we define
\begin{equation} \theta^u(z) =
-\sum_i :\be^{\rho(u)(x_i)}(z)\ga^{x'_i}(z):\ ,\end{equation} which is analogous to (1.3). An OPE computation
shows that
$$\theta^u(z)\theta^v(w) \sim B(u,v) (z-w)^{-2} + \theta^{[u,v]}(w)(z-w)^{-1}.$$ \end{proof}


\subsection{The commutant construction} There is a way to construct
interesting vertex subalgebras known as {\it commutant subalgebras}
of a given vertex algebra, which is analogous to the commutant
construction in the theory of associative algebras.

\begin{defn} Let $\vv$ be a vertex algebra and let $A$ be any subset of $\vv$.
The commutant of $A$ in $\vv$, denoted by $Com(A,\vv)$, is defined
to be the set of vertex operators $v(w)\in\vv$ which strictly
commute with the elements of $A$, that is, $[a(z),v(w)]=0$ for all
$a(z)\in A$.
\end{defn}
It follows from Lemma 2.3 that $[a(z),v(w)]=0$ iff $a(z)\ci v(z)=0$
for all $n\ge 0$, so
$$Com(A,\vv)=\{v(z)\in \vv|\ a(z)\ci v(z) = 0,\
\forall a(z)\in A,\ n\ge 0\}.$$ For any $A$, $Com(A,\vv)$ is a
vertex subalgebra, and $Com(A,\vv)= Com(\A,\vv)$, where
$\A\subset\vv$ is the vertex subalgebra generated by $A$. We regard
$\vv$ as a module over $\A$ via the left regular action, and we
regard $Com(\A,\vv)$, which will be denoted by $\vv^{\A_+}$, as the
invariant subalgebra.

If $\A$ is a homomorphic image of a current algebra
$\mathcal{O}(\G,B)$, $\A$ is generated by the subset $A = \{u(z)|\
u\in\G\}$. Hence $\vv^{\A_+} = \vv^{\G[t]}$. Consider the case $\vv = \Se(V)$ and $\A
=\Theta(\G) = \hat{\rho}\big(\mathcal{O}(\G,B)\big)$, where $\G$ is
semisimple and $V$ is a finite-dimensional $\G$-module. We claim
that generically, $\Se(V)^{\Theta_+}$ is a conformal vertex algebra.
Suppose first that $\G$ is simple, so
$$B(u,v)=-Tr\big(\rho(u)\rho(v)\big)=\lambda K(u,v)$$ for some scalar $\lambda\in
\CC$. If $\lambda\neq -\frac{1}{2}$, $\mathcal{O}(\G,\lambda K)$ has the Virasoro
element $L_{\mathcal{O}}(z)$ given by (2.13). An OPE calculation
shows that
\begin{equation}\LL(z)=L_{\Se}(z)-\hat{\rho}\big(L_{\mathcal{O}}(z)\big)\end{equation}
is a Virasoro element of central charge
$\frac{(2\lambda+2)dim(\G)}{2\lambda+1}$. In particular, if $V$ is the adjoint
module of $\G$, $\lambda=-1$ and $\LL(z)$ has central charge 0.

\begin{lem} $\LL(z)$ lies the commutant $\Se(V)^{\Theta_+}$.
Moreover, $\big(\Se(V)^{\Theta_+},\LL(z)\big)$ is a conformal vertex
algebra, and for any $a(z)\in\Se(V)^{\Theta_+}$, the OPEs of
$\LL(z)a(w)$ and $L_{\Se}(z)a(w)$ coincide.
\end{lem}
\begin{proof}
Clearly each $\theta^u(z)$ is primary of weight 1 relative to both
$L_{\Se}(z)$ and $\hat{\rho}\big(L_{\mathcal{O}}(z)\big)$. It
follows that $\LL(z)$ commutes with each $\theta^u(z)$. Hence
$\LL(z)\in\Se(V)^{\Theta_+}$.

It follows from (2.13) that any $a(z)\in\Se(V)^{\Theta_+}$ will
satisfy the OPE $\hat{\rho}\big(L_{\mathcal{O}}(z)\big)a(w)\sim 0$,
so the OPEs of $L_{\Se}(z)a(w)$ and $\LL(z)a(w)$ coincide. In
particular, the conformal weight grading on
$\big(\Se(V)^{\Theta_+},\LL(z)\big)$ coincides with the conformal
weight grading on $\Se(V)^{\Theta_+}$ inherited from the ambient
space $\big(\Se(V),L_{\Se}(z)\big)$.
\end{proof}

\begin{rem}If $\G$ is semisimple, the bilinear form
$B$ on $\G$ will be a linear combination of the Killing forms
corresponding to the various simple components of $\G$. Hence
$\hat{\rho}\big(L_{\mathcal{O}}(z)\big)$ will be a linear combination of
terms of the form (2.13) whenever it exists. Thus $\big(\Se(V)^{\Theta_+},\LL(z)\big)$ will generically be a conformal vertex algebra, although the formula for the central charge will be
more complicated.\end{rem}

Even when $\LL(z)$ is not defined, the Fourier modes $L_{\Se}(n)$
preserve $\Se(V)^{\Theta_+}$ for all $n\ge 0$. In particular, $L_{\Se}(0)$ acts by $\partial$ and $L_{\Se}(1)$ acts diagonalizably, so $\Se(V)^{\Theta_+}$ is still a {\it quasi-conformal} vertex algebra and is graded by conformal weight. 


\begin{lem} The weight zero subspace $\Se(V)^{\Theta_+}_0\subset \Se(V)^{\Theta_+}$ coincides with the classical invariant ring
$Sym(V^*)^{\G}$. In other words, $\Se(V)^{\Theta_+}$ is a \lq\lq
chiralization" of $Sym(V^*)^{\G}$.
\end{lem}
\begin{proof} Clearly $\Se(V)^{\Theta_+}_0\subset
\Se(V)^{\Theta_0}_0$, and $\Se(V)^{\Theta_0}_0 = Sym(V^*)^{\G}$. We
need to show that if $\omega\in \Se(V)_0$ is $\G$-invariant, then
$\omega$ is automatically $\Theta_+$-invariant as well. This is clear since $\theta^u(n)$ is homogeneous of
conformal weight $-n$, and the conformal weight grading on $\Se(V)$
is bounded below by 0.
\end{proof}

Next, we show that $\Se(V)^{\Theta_+}$ always contains a canonical element which is analogous to the Euler operator $\sum_i \be^{x_i}\ga^{x'_i}\in\D(V)^{\G}$.

\begin{lem} For any semisimple $\G$
and finite-dimensional module $V$, the vertex operator
$$v(z) = \sum_i :\be^{x_i}(z)\ga^{x'_i}(z):$$ lies in the
commutant $\Se(V)^{\Theta_+}$. Here $x_1,\dots,x_n$ is any basis of
$V$ and $x'_1,\dots,x'_n$ is the corresponding dual basis of $V^*$.
\end{lem}
\begin{proof}
Clearly $v(z)\in\Se(V)^{\Theta_0}$ since the pairing between $V$ and $V^*$ is $\G$-invariant. It suffices to shows that for any $u\in\G$, $\theta^u(z) \circ_1 v(z) = 0$. An OPE calculation shows that for any $u\in \G$, $\theta^u(z)\circ_1 v(z) = - Tr(\rho(u))$, which vanishes since $\G$ is semisimple. \end{proof}

\begin{rem}
The element $v(z)$ given by Lemma 2.13 satisfies the OPE
$$v(z)v(w)\sim -dim(V)(z-w)^{-2},$$ so it generates a
copy of the Heisenberg vertex algebra of central charge $-dim(V)$
inside $\Se(V)^{\Theta_+}$. \end{rem}

Suppose next that the $\G$-module $V$ admits a symmetric, $\G$-invariant bilinear form. Recall that
$\D(V)^{\G}$ has an $sl(2)$-module structure given by:
$$\psi(h) = \sum_i \be^{x_i}\ga^{x'_i},\ \
\psi(x) = \frac{1}{2} \sum_i \ga^{x'_i}\ga^{x'_i},\ \ \psi(y) =
-\frac{1}{2} \sum_i \be^{x_i}\be^{x_i},$$ where $x_1,\dots,x_n$ is an
orthonormal basis of $V$ and $x'_1,\dots,x'_n$ is the corresponding
dual basis of $V^*$.
\begin{lem}
For any semisimple Lie algebra $\G$ and any $\G$-module $V$ equipped with a symmetric, $\G$-invariant bilinear form, the homomorphism $\psi: sl(2)\rightarrow \D(V)^{\G}$ induces a
vertex algebra homomorphism
$$\hat{\psi}:\mathcal{O}\big(sl(2),-\frac{dim(V)}{8}K\big)\rightarrow
\Se(V)^{\Theta_+},$$ sending
$$h(z)\mapsto v^h(z)= \sum_i :\be^{x_i}(z)\ga^{x'_i}(z):\ ,$$
$$x(z)\mapsto v^x(z) = \frac{1}{2}\sum_i :\ga^{x'_i}(z)\ga^{x'_i}(z):\ ,$$
$$y(z)\mapsto v^y(z) = -\frac{1}{2}\sum_i :\be^{x_i}(z)\be^{x_i}(z):\ ,$$
where $K$ is the Killing form of $sl(2)$. Note that $v^h(z)$
coincides with $v(z)$ given by Lemma 2.13.
\end{lem}
\begin{proof}
This is a straightforward OPE calculation.
\end{proof}

\section{Category $\Re$}
In this section we introduce a certain category $\Re$ of vertex algebras, together with a functor from $\Re$ to the category of supercommutative rings. $\Re$ contains all vertex algebras of the form $\Se(V)$, $\E(V)$, and $\mathcal{O}(\G,B)$ and is closed under
taking subalgebras, so $\Theta$, $\Se(V)^{\Theta_+}$, $\A$, and $\Se(V)^{\A_+}$ lie in $\Re$ as well. This functor provides a bridge between vertex algebras and commutative algebra, and it allows us to answer structural question about vertex algebras $\vv\in\Re$ by using the tools of commutative algebra. 

\begin{defn}
Let $\Re$ be the category of pairs $(\vv,deg)$, where $\vv$ is a
vertex algebra equipped with a $\mathbb{Z}_{\ge 0}$-filtration
\begin{equation}\vv_{(0)}\subset\vv_{(1)}\subset\vv_{(2)}\subset \cdots,\ \ \ \vv = \bigcup_{k\ge 0}
\vv_{(k)}\end{equation} such that $\vv_{(0)} = \mathbb{C}$, and for all
$a\in \vv_{(k)}$, $b\in\vv_{(l)}$, we have
\begin{equation} a\ci b\in\vv_{(k+l)},\ \ \ for\
n<0,\end{equation}
\begin{equation} a\ci b\in\vv_{(k+l-1)},\ \ \ for\
n\ge 0.\end{equation} 
Here $\vv_{(k)}:=0$ for $k<0$.
A nonzero element $a(z)\in\vv$ is said to have
degree $d$ if $d$ is the minimal integer for which
$a(z)\in\vv_{(d)}$. Morphisms in $\Re$ are morphisms of vertex
algebras which preserve the above filtration.
\end{defn}


Filtrations on vertex algebras satisfying (3.2)-(3.3) were introduced in \cite{Li2} and are known as {\it good increasing filtrations}. If $\vv$ possesses such a filtration, it follows from (2.6)-(2.8) that the associated graded object $$gr(\vv) = \bigoplus_{k>0}\vv_{(k)}/\vv_{(k-1)}$$ is a
$\mathbb{Z}_{\ge 0}$-graded associative, supercommutative algebra with a
unit $1$ under a product induced by the Wick product on $\vv$.
Moreover, $gr(\vv)$ has a derivation $\partial$ of degree zero
(induced by the operator $\partial = \frac{d}{dz}$ on $\vv$), and
for each $a\in\vv_{(d)}$ and $n\ge 0$, the operator $a\ci$ on $\vv$
induces a derivation of degree $d-1$ on $gr(\vv)$. Finally, these
derivations give $gr(\vv)$ the structure of a {\it vertex Poisson
algebra},  ie, a graded associative,
super-commutative algebra $\A$ equipped with a derivation
$\partial$, and a family of derivations $a(n)$ for each $n\ge 0$ and
$a\in \A$ \cite{Li2}\cite{FBZ}. 

We do {\it not} require the filtration on $\vv$ to come from
a $\mathbb{Z}_{\ge 0}$-grading
$$\vv=\bigoplus_{k\ge 0} V^{(k)}$$ where $\vv_{(k)}=\oplus_{i=0}^k
V^{(i)}$. If $\vv$ does possess such a grading, we will say that
$\vv$ is {\it graded by degree}. If $\A$ is a vertex subalgebra of
$\vv$, the filtration (3.1) on $\vv$ induces a filtration
$$\A_{(0)}\subset \A_{(1)}\subset \A_{(2)}\subset \cdots$$ on $\A$,
where $\A_{(k)}=\A\bigcap \vv_{(k)}$. With respect to this
filtration, $(\A,deg)$ lies in $\Re$. In general, if $\vv$ is graded
by degree, a subalgebra $\A$ need not be graded by degree.

In general, there is no natural linear map from $\vv\rightarrow
gr(\vv)$, but we do have projections
\begin{equation}\phi_d:\vv_{(d)}\rightarrow \vv_{(d)}/\vv_{(d-1)}\subset
gr(\vv)\end{equation} for $d\ge1$. If $a,b\in gr(\vv)$ are
homogeneous of degrees $r,s$ respectively, and $a(z)\in\vv_{(r)}$,
$b(z)\in\vv_{(s)}$ are vertex operators such that
$\phi_r\big(a(z)\big)=a$ and $\phi_s\big(b(z)\big)=b$, it follows
that $\phi_{r+s}\big(:a(z)b(z):\big) = ab$.

Let $\mathcal{R}$ denote the category of $\mathbb{Z}_{\ge 0}$-graded supercommutative rings equipped with a derivation $\partial$ of degree 0, which we shall call $\partial$-rings.

\begin{lem} If $(\vv,deg)$ and $(\vv',deg')$ lie in $\Re$, and
$f:\vv\rightarrow \vv'$ is a morphism in $\Re$, $f$ induces a
homomorphism of $\partial$-rings $gr(f):gr(\vv)\rightarrow
gr(\vv')$. In particular, the assignment $(\vv,deg)\mapsto gr(\vv)$
is a functor from $\Re$ to $\mathcal{R}$.
\end{lem}

\begin{proof} Let $a\in gr(\vv)$ be homogeneous of degree $r$, and
let $a(z)\in\vv$ be any vertex operator of degree $r$ such that
$\phi^{\vv}_r\big(a(z)\big)=a$. We define $$gr(f)(a) = \phi^{\vv'}_r
f\big(a(z)\big).$$ If $a'(z)\in\vv_{(r)}$ is another vertex operator
such that $\phi^{\vv}_{r}\big(a'(z)\big)=a$, then
$deg\big(a(z)-a'(z)\big)<r$. Since $f$ preserves degree,
$f\big(a(z)-a'(z)\big)$ also has degree less than $r$, so
$\phi^{\vv'}_r f\big(a(z)\big)-\phi^{\vv'}_r f\big(a'(z)\big)=0$.
Hence $gr(f)$ is well-defined. To see that $gr(f)$ is a ring
homomorphism, let $b\in\vv$ be homogeneous of degree $s$ and choose
$b(z)\in\vv_{(s)}$ such that $\phi^{\vv}_s\big(b(z)\big)=b$. Then
$:a(z)b(z):\in\vv_{(r+s)}$ and
$\phi^{\vv}_{(r+s)}\big(:a(z)b(z):\big)=ab$. It follows that
$$gr(f)(ab) = \phi^{\vv'}_r f\big(:a(z)b(z):\big) = \phi^{\vv'}_r
\big(:f\big(a(z)\big)f\big(b(z)\big):\big)$$ $$=\phi^{\vv'}_r
f\big(a(z)\big)\phi^{\vv'}_r f\big(a(z)\big)= gr(f)(a) gr(f)(b).$$
The same argument shows that
$gr(f)\big(\partial(a)\big)=\partial\big(gr(f)(a)\big)$. Checking
that $gr$ respects compositions of mappings and that $gr(id_{\vv}) =
id_{gr(\vv)}$ is also straightforward.
\end{proof}

A vertex algebra $\vv$ is said to be {\it freely generated} by an
ordered collection $\{a_i(z)|\ i\in I\}$, if the collection of iterated Wick products $$\{:a_{i_1}(z)\cdots
a_{i_r}(z):|\ i_1\le\cdots \le i_r\}$$ forms a {\it basis} of $\vv$ \cite{KdS}. For example, if we fix a basis $u_1,\dots,u_n$ for a Lie algebra $\G$, $\mathcal{O}(\G,B)$ is freely generated by the collection $$\{\partial^k u_i(z)|\ i=1,\dots,n,\ k\ge 0\},$$ which we order by declaring $\partial^k u_i(z)>\partial^l u_j(z)$ if $i>j$, or $i=j$ and $k>l$. 

Similarly, if $x_1,\dots,x_n$ is a basis for a vector space $V$ and $x'_1,\dots,x'_n$ is the corresponding dual basis for $V^*$, $\Se(V)$ is freely generated by the collection
$$\{\partial^k\be^{x_i},\partial^k\gamma^{x'_i}(z)|\ i=1,\dots,n,\ k\ge 0\},$$ which we order in the obvious way.

\begin{lem} Suppose that $\vv$ is freely generated by an
ordered collection of vertex operators $\{a_i(z)|\ i\in I\}$. Assign
each $a_i(z)$ degree $d_i>0$, and define a linear $\mathbb{Z}_{\ge
0}$-grading $\vv=\bigoplus_{n\ge 0} V^{(n)}$ by declaring
\begin{enumerate}
\item $\vv^{(0)}=\mathbb{C}$
\item $\vv^{(n)}$ is spanned by the collection of vertex operators
$:a_{i_1}(z)\cdots a_{i_r}(z):$ for which $d_{i_1} +\cdots +d_{i_r} = n$.
\end{enumerate}
If the generators $a_i(z)$ satisfy (3.2)-(3.3), then $(\vv,deg)$ lies in $\Re$.
\end{lem}

\begin{proof} This is a special case of Theorem 4.6 of \cite{Li2}. A straightforward induction on the degree and on the number of derivatives
shows that any pair of homogeneous vertex operators
$a(z),b(z)\in\vv$ of degrees $d_a,d_b$, respectively, satisfies
$a(z)\ci b(z)\in\vv_{d_a+d_b-1}$ for $n\ge 0$. Hence (3.3)
holds for all of $\vv$. (3.2) then follows immediately from
(2.6).
\end{proof}

\begin{cor} If we equip $\mathcal{O}(\G,B)$ with the grading
$deg(\partial^k u(z))=1$ for all $u\in\G$ and $k\ge 0$,
$\mathcal{O}(\G,B)$ lies in $\Re$ and is graded by degree.
$\mathcal{O}(\G,B)^{(n)}$ is spanned by the collection of vertex
operators
$$:\partial^{k_1}u_1(z)\cdots\partial^{k_n} u_n(z):\ ,$$ where $u_i\in\G$ and $k_i\ge 0$.
Moreover, $gr\big(\mathcal{O}(\G,B)\big)$ is the polynomial algebra
$$Sym \big(\bigoplus_{k\ge 0} \G_k\big),$$ where $\G_k$ is the
copy of $\G$ spanned by the collection of vertex operators
$\{\partial^k u(z)\}$ which are linear in $u\in\G$.
\end{cor}

\begin{cor} If we equip $\Se(V)$ with the grading $$deg(\partial^k \Gx(z))=deg(\partial^k\Bx(z))
=1$$ for all $x\in V$, $x'\in V^*$, and $k\ge 0$, $\Se(V)$ lies in
$\Re$ and is graded by degree. $\Se(V)^{(n)}$ is spanned by the
collection
$$:\partial^{k_1}\be^{x_1}(z)\cdots\partial^{k_r}\be^{x_r}(z)\partial^{l_1}\ga^{x'_1}(z)\cdots
\partial^{l_s}\ga^{x'_s}(z):\ ,$$ and $gr(\Se(V))$ is the polynomial algebra $$Sym
\big(\bigoplus_{k\ge 0} (V_k\oplus V^*_k)\big),$$ where $V_k$ and
$V^*_k$ are the copies of $V$ and $V^*$, spanned by the collections
$\{\partial^k\Bx(z)\}$ and $\{\partial^k \Gx(z)\}$, which are linear
in $x\in V$ and $x'\in V^*$, respectively.
\end{cor}

If $(\vv,deg)\in\Re$, we may rescale $deg$ by a factor of $m$ for
any positive integer $m$, and the resulting pair $(\vv,m\cdot deg)$
will still lie in $\Re$. When $V$ is a $\G$-module via
$\rho:\G\rightarrow End(V)$, the map
$$\hat{\rho}:\mathcal{O}(\G,B)\rightarrow \Se(V)$$
given by Lemma 2.8 is a morphism in $\Re$ if we double the above
grading on $\mathcal{O}(\G,B)$, so that $deg\big(u(z)\big) = 2$ for
all $u\in\G$. Likewise, if $V$ admits a symmetric, $\G$-invariant bilinear form, the map
$$\hat{\psi}:\mathcal{O}\big(sl(2),-\frac{dim
(V)}{8}K\big)\rightarrow \Se(V)^{\Theta_+}$$ given by Lemma 2.15 is
a morphism in $\Re$ if the degree grading on
$\mathcal{O}\big(sl(2),-\frac{dim(V)}{8}K\big)$ is doubled.

Let $V$ be a $\G$-module as above, and let $\Se(V)^{\Theta_0}\subset \Se(V)$ denote the $\G$-invariant subalgebra which is annihilated by $\{\theta^u(0)|\ u\in\G\}$. Clearly $\Se(V)^{\Theta_+}\subset \Se(V)^{\Theta_0}$, and these vertex algebras both lie in $\Re$ as subalgebras of $\Se(V)$ with the induced filtration. $\Se(V)^{\Theta_0}$ is graded by degree as a subalgebra of $\Se(V)$
since each $\{\theta^u(0)|\ u\in\G\}$ is a homogeneous derivation of degree 0 on $\Se(V)$. The associated graded algebra $gr\big(\Se(V)^{\Theta_0}\big)$ is the classical invariant
ring $$Sym\big(\bigoplus_{k\ge 0} (V_k\oplus V^*_k)\big)^{\G}.$$
However, for $n>0$, $\theta^u(n)$ decomposes into
homogeneous components of degrees 0 and $-2$, so $\Se(V)^{\Theta_+}$
is {\it not} graded by degree as a subalgebra of $\Se(V)$. 

Similarly, if $V$ admits a symmetric, $\G$-invariant bilinear form, so that $\Se(V)^{\Theta_+}$ contains the subalgebra $\A = \hat{\psi}\big(\mathcal{O}\big(sl(2),-\frac{dimV}{8}K\big)\big)$, the $sl(2)$-invariant space $\Se(V)^{\A_0}$ is graded by degree but $\Se(V)^{\A_+}$ is not, since the operators $v^u(n)$ decompose into homogeneous components of degrees $0$ and $-2$, for $n>0$.

The key feature of $\Re$ is that vertex algebras $\vv\in\Re$ have
the following {\it reconstruction property}: we can write down a set
of strong generators for $\vv$ as a vertex algebra just by knowing
the ring structure of $gr(\vv)$. We say that the collection $\{a_i|\
i\in I\}$ generates $gr(\vv)$ as a $\partial$-ring if the collection
$\{\partial^k a_i|\ i\in I,\ k\ge 0\}$ generates $gr(\vv)$ as a
graded ring.

\begin{lem} Let $\vv$ be a vertex algebra in $\Re$. Suppose that
$gr(\vv)$ is generated as a $\partial$-ring by a collection
$\{a_{i}|\ i\in I\}$, where $a_i$ is homogeneous of degree $d_i$.
Choose vertex operators $a_i(z)\in \vv_{(d_i)}$ such that
$\phi_{d_i}\big(a_i(z)\big) = a_i$. Then $\vv$ is strongly generated
by the collection $\{a_i(z)|\ i\in I\}$.
\end{lem}
\begin{proof}
Let $\vv'\subset\vv$ denote the linear subspace spanned by the
monomials
$$:\partial^{k_1}a_{i_1}(z)\cdots\partial^{k_r}a_{i_r}(z):\ .$$ We
need to prove that $\vv'=\vv$; we proceed by induction on degree.
The statement is trivial in degree 0, so assume it for degree less
than $d$. Let $\omega(z)\in \vv_{(d)}$ and let $\omega\in gr(\vv)$
denote the image of $\omega(z)$ under $\phi_d:\vv_{(d)}\rightarrow
\vv_{(d)}/\vv_{(d-1)}$. Since $gr(\vv)$ is generated as a
$\partial$-ring by $a_i$, we can write
$$\omega = \sum_{K,I} \lambda_{K,I}\partial^{k_1}a_{i_1}\cdots\partial^{k_r}a_{i_r},$$
where the sum is over all monomials in $gr(\vv)$ for which
$d_1+\dots+ d_r = d$. Let
$$\omega'(z) = \sum_{K,I}
\lambda_{K,I}:\partial^{k_1}a_{i_1}(z)\cdots\partial^{k_r}a_{i_r}(z):\
.$$ It is easy to see that $\phi_d
\big(\omega'(z)\big)=\omega=\phi_d\big(\omega(z)\big)$, so that
$\omega''(z)=\omega(z)-\omega'(z)$ has degree less than $d$. Since
$\omega'(z)\in \vv'$, we have $\omega''(z)\equiv\omega(z)$ modulo
$\vv'$. The claim follows by induction.\end{proof}


Recall that for all $a(z)\in\vv_{(d)}$ and $n\ge 0$, $a(z)\ci$
induces a derivation of degree $d-1$ on $gr(\vv)$, and these maps
give $gr(\vv)$ the structure of a vertex Poisson algebra \cite{Li2}{\cite{FBZ}.
However, this structure may be trivial in the sense that all these
maps on $gr(\vv)$ are zero. If $\vv$ is abelian (ie, $[a(z),b(w)]=0$
for all $a,b\in\vv$) this will always be the case, but it may be
true even if $\vv$ is nonabelian. However, we will modify this
construction slightly to obtain a {\it non-trivial} vertex Poisson
algebra structure on $gr(\vv)$ whenever $\vv$ is not abelian. Define
$$
k= k(\vv,deg)=sup\{j\geq1\vert
\vv_{(r)}\ci\vv_{(s)}\subset \vv_{(r+s-j)}~\forall r,s,n\geq0\}.
$$ 
It follows easily that $k$ is finite iff $\vv$ is not abelian.


\begin{lem}
Let $(\vv,deg)\in\Re$ and let $k=k(\vv,deg)$ be as above. For each
$a(z)\in\vv$ of degree $d$ and $n\ge 0$, the operator $a(z)\ci$ on
$\vv$ induces a homogeneous derivation $a(n)_{Der}$ on $gr(\vv)$ of
degree $d-k$, defined on homogeneous elements $b$ of degree $r$ by
\begin{equation} a(n)_{Der}(b) = \phi_{r+d-k}\big(a(z)\ci
b(z)\big).\end{equation} Here $b(z)\in\vv$ is any vertex operator of
degree $r$ such that $\phi_r \big(b(z)\big) = b$.
\end{lem}
\begin{proof}
If $b'(z)\in \vv$ is another vertex operator of degree $r$ such that
$\phi_r\big(b'(z)\big)=b$, then $deg\big((b(z)-b'(z)\big)<r$. Using
(3.5), it follows that $$\phi_{r+d-k}\big(a(z)\ci
b'(z)\big)-\phi_{r+d-k}\big(a(z)\ci b(z)\big)=
\phi_{r+d-k}\big(a(z)\ci\big(b(z)-b'(z)\big)\big)=0.$$ Hence
$a(n)_{Der}$ is well-defined, and is clearly homogeneous of degree
$d-k$. It remains to show that for any homogeneous $b,c\in\vv$ of
degrees $r,s$ respectively, we have
$$a(n)_{Der}(bc)-\big(a(n)_{Der}(b)\big)c-b
\big(a(n)_{Der}(c)\big)=0.$$ Let $b(z),c(z)\in\vv$ be vertex
operators of degrees $r,s$ respectively, such that
$\phi_r\big(b(z)\big)=b$ and $\phi_s\big(c(z)\big)=c$, so that
$\phi_{r+s}\big(:b(z)c(z):\big)=bc$. Hence $a(n)_{Der}(bc) =
\phi_{r+s+d-k} a(z)\ci \big(:b(z)c(z):\big)$. Similarly,
$$\big(a(n)_{Der}(b)\big)c = \phi_{r+s+d-k}\big( :\big(a(z)\ci
b(z)\big) c(z):\big),$$ and
$$b
\big(a(n)_{Der}(c)\big) = \phi_{r+s+d-k}\big( :b(z)\big(a(z)\ci
c(z)\big):\big).$$ Hence $a(n)_{Der}(bc)-\big(a(n)_{Der}(b)\big)c-b
\big(a(n)_{Der}(c)\big)$ is equal to $$\phi_{r+s+d-k}\bigg(a(z)\ci
\big(:b(z)c(z):\big)-\ :\big(a(z)\ci b(z)\big) c(z):\ -\ :b(z)
\big(a(z)\ci c(z)\big):\bigg).$$ Using (2.7), we see that this
expression is equal to
$$\phi_{r+s+d-k}\bigg(\sum_{i=1}^n \binom ni
\big(a(z)\circ_{n-i}b(z)\big)\circ_{i-1}c(z)\bigg).$$ By (3.5),
$a(z)\circ_{n-i}b(z)\in \vv_{r+d-k}$, and
$\big(a(z)\circ_{n-i}b(z)\big)\circ_{i-1}c(z) \in \vv_{r+s+d-2k}$ by
applying (3.5) again. Since $k\ge 1$, the claim follows.
\end{proof}

Clearly the maps $\{a(n)_{Der}|\ a\in\vv,\ n\ge 0\}$ give $gr(\vv)$
the structure of a vertex Poisson algebra, and this structure is
nontrivial whenever $k$ is finite. Note that if $(\vv,deg)$ lies in
$\Re$ and we rescale the degree by a factor of $m$, $k(\vv,m\cdot
deg)= m\cdot k(\vv,deg)$. If $\vv$ is strongly generated by a set
$\{a_i(z)|\ i\in I\}$ of vertex operators of degrees $d_i$
satisfying the conditions of Lemma 3.3, it is easy to see that
$k(\vv,deg)$ is the minimum value of
$$deg (a_i(z))+deg(a_j(z))-deg \big(a_i(z)\ci a_j(z)\big),$$ where $i,j$
range over $I$ and $n\ge 0$. It follows from the OPE formulas (2.11)
and (2.14) that
\begin{equation} k(\mathcal{O}(\G,B),deg)=1,\ \ \ k(\Se(V),deg)=2.\end{equation}

\begin{lem} Let $(\vv,deg)\in \Re$, and suppose $a(z),b(z)\in \vv$
are vertex operators of degrees $r$ and $s$ such that
$[a(z),b(z)]=0$. Then for all $n,m\ge 0$, $a(n)_{Der}$ and
$b(m)_{Der}$ commute as operators on $gr(\vv)$.
\end{lem}
\begin{proof} Let $c\in gr(\vv)$ be homogeneous of degree $t$, and let
$c(z)\in\vv$ be a vertex operator of degree $t$ such that
$\phi_t\big(c(z)\big)=c$. Then $b(m)_{Der}(c) = \phi_{s+t-k}
\big(b(z)\circ_m c(z)\big)$. Likewise,
$$a(n)_{Der}\big(b(m)_{Der}(c)\big)=\phi_{r+s+t-2k} \big(a(z)\ci
\omega(z)\big),$$ where $\omega(z)$ is any vertex operator of degree
$s+t-k$ such that $\phi_{s+t-k}\big(\omega(z)\big)=b(m)_{Der}(c)$.
We may take $\omega(z)=b(z)\circ_m c(z)$. Then
$$a(n)_{Der}\big(b(m)_{Der}(c)\big)=\phi_{r+s+t-2k} \big(a(z)\ci
\big(b(z)\circ_m c(z)\big)\big),$$ and similarly,
$$b(m)_{Der}\big(a(n)_{Der}(c)\big)=\phi_{r+s+t-2k}
\big(b(z)\circ_m \big(a(z)\ci c(z)\big)\big).$$ It follows that
\begin{equation}[a(n)_{Der},b(m)_{Der}](c)=\phi_{r+s+t-2k}\big([a(n),b(m)]\big(c(z)\big)\big).\end{equation}
Since $a(z),b(z)$ commute, it follows that $[a(n),b(m)]=0$ for all
$n,m\ge 0$, which proves the claim.
\end{proof}

\subsection{Commutants in $\Re$}
Let $(\vv,deg)\in\Re$, $k=k(\vv,deg)$, and let $\A$ be a subalgebra
of $\vv$ which is a homomorphic image of a current algebra
$\mathcal{O}(\G,B)$. We would like to use the filtration $deg$ on
$\vv$ to study the commutant $\vv^{\A_+}$. Suppose that for each
$u\in\G$, $u(z)\in \A$ has degree $k$, so that the derivations
$\{u(n)_{Der}|\ n\ge 0\}$ on $gr(\vv)$ are homogeneous of degree 0
by Lemma 3.7.

\begin{lem} The derivations $\{u(n)_{Der}|\ n\ge 0\}$ form a
representation of $\G[t]$ on $gr(\vv)$. Moreover, the actions of
$\G[t]$ on $\vv$ and $gr(\vv)$ are compatible in the sense that for
any $\omega(z)\in\vv$ of degree $r$, we have
\begin{equation} u(n)_{Der}\phi_r \big(\omega(z)\big) =
\phi_r\circ u(n)\big(\omega(z)\big).\end{equation}
\end{lem}
\begin{proof} Let $\omega\in gr(\vv)\big)$ be homogeneous of degree $r$,
and let $\omega(z)\in\vv$ be a vertex operator of degree $r$ such
that $\phi_r\big(\omega(z)\big)=\omega$. Using (3.7) and the fact
that each $u(n)_{Der}$ has degree 0, we have
$$[u(n)_{Der},v(m)_{Der}](\omega) =
\phi_{r}\big([u(n),v(m)]\big(\omega(z)\big)\big)$$
$$=\phi_{r}\big([u,v](n+m)\big(\omega(z)\big)\big)=[u,v](n+m)_{Der}(\omega).$$
\end{proof}

Since each $u(n)_{Der}$ is degree-homogeneous, the invariant space
$gr(\vv)^{\A_+}$ under this action is graded by degree as a subalgebra of
$gr(\vv)$. Moreover, $gr(\vv)^{\A_+}$ is closed under
$\partial$ since $\vv^{\A_+}$ is a vertex algebra and $\partial$ is
homogeneous of degree 0. By functoriality, the inclusion of vertex
algebras $\vv^{\A_+}\subset \vv$ gives rise to an injective ring
homomorphism $gr(\vv^{\A_+})\hookrightarrow gr(\vv)$ whose image
clearly lies in $gr(\vv)^{\A_+}$. Hence we have a canonical
injection
\begin{equation} \Gamma: gr(\vv^{\A_+})\hookrightarrow
gr(\vv)^{\A_+}\end{equation} which is a homomorphism of
$\partial$-rings.

\subsection{A strategy for computing $\vv^{\A_+}$}
Let $R\subset gr(\vv)^{\A_+}$ denote the image of $gr(\vv^{\A_+})$ under $\Gamma$. The problem of finding a set of generators $\{a_i|\ i\in I\}$ for $R$ as a $\partial$-ring is a problem in commutative algebra. Solving this problem allows us find generators for the vertex algebra $\vv^{\A_+}$ as well.
Since $\Gamma$ maps $gr(\vv^{\A_+})$ isomorphically onto $R$, these generators correspond to elements of
$gr(\vv^{\A_+})$, which in turn come from vertex operators
$\{a_i(z)|\ i\in I\}$ in $\vv^{\A_+}$ such that
$\phi_{d_{i}}\big(a_i(z)\big)=a_i$. Here $d_i = deg(a_i)$. By Lemma 3.6, this collection strongly
generates $\vv^{\A_+}$. {\it In particular, $\vv^{\A_+}$ is (strongly) finitely generated as a vertex algebra whenever $R$ is finitely generated as a $\partial$-ring.} 

In our main example, we will find a finite set of generators for $gr(\vv)^{\A_+}$ as a $\partial$-ring. These generators correspond to vertex operators in $\vv^{\A_+}$, so in this case $\Gamma$ is surjective (and hence an isomorphism), and we obtain a finite set of generators for $\vv^{\A_+}$ as well.


Consider the case $\vv=\Se(V)$ and $\A=\Theta(\G)$, where $\G$ is semisimple and $V$ is a
finite-dimensional $\G$-module. In this case, $deg(\theta^u(z)) = 2
= k$, so each $\theta^u(n)_{Der}$ is homogeneous of degree 0 and
$gr(\Se(V))$ is a $\G[t]$-module by Lemma 3.9. For notational
simplicity, we denote $gr\big(\Se(V)\big)$ by $P$, and we denote the images of $\partial^k \Bx(z),\partial^k
\Gx(z)$ in $P$ by
$\Bx_k$ and $\Gx_k$, respectively. The action of $\theta^u(n)_{Der}$
on the generators of $P$ is given by
\begin{equation}\theta^u(n)_{Der}(\Bx_k) = c^n_k
\be^{\rho(u)(x)}_{k-n},\ \ \ \ \ \theta^u(n)_{Der}(\Gx_k) = c^n_k
\ga^{\rho^*(u)(x')}_{k-n},\end{equation} where $c^n_k =
k(k-1)\cdots(k-n+1)$, for $n,k\ge 0$. Clearly $c^0_k = 1$ and
$c^n_k=0$ for $n>k$.

If $V$ admits a symmetric, $\G$-invariant bilinear form, so that by Lemma 2.15, $\Se(V)^{\Theta_+}$ contains the subalgebra $\A = \hat{\psi}\big(\mathcal{O}\big(sl(2),-\frac{dim(V)}{8}K\big)\big)$, the operators $\{v^u(k)_{Der}|\ u=x,y,h,\ k\ge
0\}$ on $P$ form a representation of the Lie algebra $sl(2)[t]$ by
derivations of degree 0. In terms of an orthonormal basis of $V$,
the action is given by:

\begin{equation}v^h(n)_{Der}(\be^{x_i}_k) = -c^n_k \be^{x_i}_{k-n},\ \
v^h(n)_{Der}(\ga^{x'_i}_k) = c^n_k \ga^{x'_i}_{k-n},\end{equation}
\begin{equation}v^x(n)_{Der}(\be^{x_i}_k) = -\frac{1}{2}c^n_k \ga^{x'_i}_{k-n},\ \
v^x(n)_{Der}(\ga^{x'_i}_k) = 0,\end{equation}
\begin{equation}v^y(n)_{Der}(\be^{x_i}_k) = 0,\ \
v^y(n)_{Der}(\ga^{x'_i}_k) = -\frac{1}{2}c^n_k
\be^{x'_i}_{k-n}.\end{equation} 

We denote the invariant space $gr\big(\Se(V)\big)^{\A_+}$ by $P^{\A_+}$. Our main task is to describe $P^{\A_+}$ as a $\partial$-ring in the case where $\G= sl(2)$ and $V$ is the adjoint module. For this purpose, it is useful to define another $\mathbb{Z}_{\ge 0}$-grading on $P$
which we call {\it level}; each $\Bx_k$ and $\Gx_k$ has level $k$.
It is clear from (3.11-3.13) that each $v^u(n)_{Der}$ is homogeneous
of level $-n$, so $P^{\A_+}$ is graded by level. 
In addition to the gradings $deg$ and $lev$, $P$ has various auxiliary $\mathbb{Z}_{\ge 0}$-gradings which will be useful. An essential argument is to show that the condition $\omega\in P^{\A_+}$ implies that the projection of $\omega$ onto certain homogeneous subspaces is nonzero (see Lemmas 4.15-4.17).

\subsection{Grobner bases}
$P^{\A_+}$ has an additional feature; it is a subalgebra of the classical invariant ring $P^{\A_0} = P^{sl(2)}$. In the case $\G = sl(2) = V$, $P^{\A_0}$ can be exhibited as a quotient $F/I$,
where $F$ is a polynomial algebra on countably many variables, and
$I$ is a countably generated ideal. Regarding
$P^{\A_+}$ as a subalgebra of $F/I$, we can study it using the tools of
commutative algebra. In particular, we can find a Grobner basis for
$I$ and a corresponding normal form for elements of $P^{\A_0}$.
By passing back and forth between the description of $P^{\A_+}$ as
a subalgebra of $P$ and as a subalgebra of $F/I$, we will give a
complete description of $P^{\A_+}$.

Even though $P^{\A_0}$ is not finitely generated, it
has a natural filtration by finitely generated subalgebras. $P$ is
filtered by the subalgebras
$$P_N = Sym \big(\bigoplus_{k=0}^N (V_k\oplus V^*_k)\big),\ \ N\ge 0,$$
which are generated by $\Bx_k,\Gx_k$ for $k=0,\dots,N$. By (3.11-3.13),
the action of $sl(2)[t]$ on $P$ preserves each $P_N$, so $P^{\A_0}$
and $P^{\A_+}$ are filtered by the subalgebras $P^{\A_0}_N =
P^{\A_0}\bigcap P_N$ and $P^{\A_+}_N = P^{\A_+}\bigcap
P_N$, respectively. Hence when working in $P^{\A_0}$ and
$P^{\A_+}$, we may always assume that we are working inside some
$P_N$ for $N$ sufficiently large. Then $P^{\A_0}_N$ will be a
quotient $F/I$ of a finitely generated polynomial ring $F$,
and we can apply the standard techniques of commutative algebra (localization, Grobner basis theory, etc) without difficulty.

We recall the definition and basic properties of Grobner bases,
following \cite{CLO}. Let $F$ be the polynomial ring $\mathbb
C[x_1,...,x_n]$, and let $I\subset F$ be an ideal. By the Hilbert
Basis Theorem, $I$ is finitely generated; let $I=\langle
g_1,..,g_k\rangle$. Fix a monomial ordering on $F$. We will always
assume that this ordering comes from ordering the generators
$$x_1<x_2<\cdots<x_n,$$ and then ordering monomials in $F$
lexicographically. For any polynomial $f\in F$, we denote the leading
term of $f$ with respect to this ordering by $lt(f)$. Let $\langle
lt(I)\rangle$ denote the monomial ideal generated by the collection
$\{lt(f)|\ f\in I\}$.

\begin{defn} We say that the collection $\{g_1,...,g_k\}$ forms a Grobner
basis for $I$ if $\langle lt(g_1),\dots,lt(g_k)\rangle = \langle
lt(I)\rangle$.
\end{defn}
The key property of Grobner bases is the following.
\begin{thm} Let $B = \{g_1,\dots g_k\}$ be a Grobner basis for $I$. Then for any $f\in F$, there is a unique $r\in F$ with the following properties: (i) No monomial appearing in $r$
is divisible by $lt(g_1),\dots, lt(g_k)$. (ii) There is a unique
$g\in I$ such that $f=g+r$.
\end{thm}
The polynomial $r$ is called the {\it normal form} of $f$. By
uniqueness, $r$ is the remainder of $f$ upon long division by the
generators $g_1,\dots,g_k$ in any order. Clearly normal forms behave
well under addition; if $r_1$ and $r_2$ are the normal forms of
$f_1$ and $f_2$, respectively, then $r_1+r_2$ is the normal form of
$f_1+f_2$. There is a procedure, known as {\it Buchberger's
algorithm}, for extending a given set of generators
$B=\{g_1,\dots,g_k\}$ for $I$ to a Grobner basis. We sketch this procedure, following the notation in \cite{CLO}.

For any two polynomials $f,g\in F$, define the $S$-polynomial
$$ S(f,g) = \frac{LCM(lt(f),lt(g))}{lt(f)} f - \frac{LCM(lt(f),lt(g))}{lt(g)} g.$$
For any two elements $g_i,g_j\in B$, let $\overline{S(g_i,g_j)}^B$ denote the remainder of $S(g_i,g_j)$ upon long division by $g_1,\dots,g_k$ (in that order). Buchberger showed that $B$ is a Grobner basis for $I$ if and only if $\overline{S(g_i,g_j)}^B = 0$ for every pair $i,j$. If $B$ is not a Grobner basis for $I$, we may adjoin all the (nonzero) polynomials of the form $\overline{S(g_i,g_j)}^B$ to the set $B$, obtaining a bigger set $B'$. This algorithm terminates after a finite number of steps, and the resulting set will be a Grobner basis for $I$.

\begin{lem}
Let $F'$ be the subalgebra $\mathbb{C}[x_1,\dots,x_s]\subset F$ for $s<n$. Suppose that $F'\bigcap I$ is empty, so the images
$\bar{x}_1,\dots,\bar{x}_s$ of $x_1,\dots,x_s$ in $F/I$ are
algebraically independent. Order the generators $$x_1<\cdots <x_s<x_{s+1}<\cdots <x_n,$$ and then
order monomials in $F$ using the standard lexicographic ordering. Choose a corresponding Grobner basis $B =\{g_1,\dots g_k\}$ for $I$. Then every element of $F'$ is already in normal form with respect to this Grobner basis.
\end{lem}

\begin{proof}
Since $F'\bigcap I$ is empty, each $g_i$ must contain monomials
which do not lie in $F'$. But any such monomial is greater than any
monomial in $F'$ in the above ordering, so the leading term of $g_i$
cannot lie in $F'$.
\end{proof}

Let $A = \mathbb{C}[y_1,\dots,y_m]$, and let $f:F\rightarrow A$ be a ring homomorphism with kernel $I$. As above, let $F' = \mathbb{C}[x_1,\dots,x_s]\subset F$, and suppose that $F'\bigcap I$ is empty. Order monomials in $F$ as in Lemma 3.12, and choose a corresponding Grobner basis $B$ for $I$, so that every element of $F'$ is in normal form. Let $A'$ and $A''$ denote the subalgebras $f(F)$ and $f(F')$ of $A$, respectively. For any $\omega\in A'$, let $\hat{\omega}\in F$ denote the normal form of the corresponding element $f^{-1}(\omega)\in F/I$. 

\begin{lem} Let $A'''$ be a subalgebra of $A$ satisfying $A''\subset A'''\subset A'$. Suppose that $A'''$ has the following property: for any $\omega\in A'''$ of positive degree, the normal form $\hat{\omega}\in F$ contains a monomial in $F'$ with nonzero coefficient. Then $A''' = A''$.
\end{lem}

\begin{proof} For any $\omega\in A'''$, define the {\it length} of $\omega$,
denoted by $l(\omega)$, to be the number of distinct monomials
appearing in $\hat{\omega}$ with nonzero coefficients. We proceed by
induction on length. If $l(\omega)=1$, $\hat{\omega}$ consists of a
single monomial, which lies in $F'$ by hypothesis. Since $F'$
corresponds isomorphically to $A''$ under $f$, we have $\omega\in A''$. Suppose
that $l(\omega)=n$. Consider the normal form $\hat{\omega}$, and let
$\mu$ be a monomial appearing in $\hat{\omega}$ which lies in $F'$.
Let $\omega' = \omega - f(\mu)$. Since the elements of $F'$ are
already in normal form and normal forms are additive, $\hat{\omega'}
= \hat{\omega}-\mu$, so $\omega'$ has length $n-1$. Since
$\omega'\equiv \omega$ modulo $A''$, the claim follows by
induction.
\end{proof}

\begin{rem} Consider the case $$A = P = Sym
\big(\bigoplus_{k\ge 0} (V_k\oplus V^*_k)\big),\ \ \
A'=P^{\A_0},\ \ \ A''' = P^{\A_+},\ \ \ A'' = P_{\tau}.$$ Here $P_{\tau}$ is a certain candidate for $P^{\A_+}$ (to be defined in the next section) which is contained in $P^{A_+}$ and
is generated by algebraically independent elements. We will use Lemma 3.13 to show that $P_{\tau}=P^{\A_+}$. \end{rem}

\section{The computation of $\Se(V)^{\A_+}$}
In this section, we prove Theorem 1.3. For $\G=sl(2)=V$, we denote
$\Se(V)$ and
$\Theta(\G)=\hat{\rho}\big(\mathcal{O}\big(sl(2),-K)\big)$ by $\Se$
and $\Theta$, respectively, and we work in the basis $x,y,h$ with
the commutation relations
\begin{equation}[x,y]=h,\ \ \ \  [h,x]=2x,\ \ \ \ [h,y]=-2y.\end{equation}
$\Theta$ is generated by the vertex operators
$$\tX(z) = 2:\Bx(z)\Gh(z):\ -\ :\Bh(z)\Gy(z):\ ,$$
$$\tY(z) = -2:\By(z)\Gh(z):\ +\ :\Bh(z)\Gx(z):\ ,$$ $$\tH(z) =
-2:\Bx(z)\Gx(z):\ +\ 2:\By(z)\Gy(z):\ .$$

In the above basis, the subalgebra $\A =
\hat{\psi}\big(\mathcal{O}(sl(2),-\frac{3}{8}K)\big)$ of
$\Se^{\Theta_+}$ is generated by the vertex operators
$$v^x(z) = \frac{1}{2}\big(:\Gh(z)\Gh(z):\ +\
:\Gx(z)\Gy(z):\big)$$
$$v^y(z) = -\frac{1}{2}\big(:\Bh(z)\Bh(z):\ +\ 4:\Bx(z)\By(z):\big)$$
$$v^h(z) = \ :\Bx(z)\Gx(z):\ +\ :\By(z)\Gy(z):\ +\ :\Bh(z)\Gh(z):\ .$$

In terms of the above basis, $P=gr(\Se)$ is the polynomial algebra generated by
$\{\be^u_k,\ga^{u'}_k|\ u=x,y,h,\ k\ge 0\}$. To simplify
notation, we drop the subscript $Der$ and denote the operators
$\theta^u(n)_{Der}$ and $v^u(n)_{Der}$ on $P$ by $\theta^u(n)$ and
$v^u(n)$, respectively. 
Define the polynomials
\begin{equation} \tau^x_0 = \phi_2\big(\theta^x(z)\big) = 2\Bx_0\Gh_0 -\Bh_0\Gy_0, \end{equation} 
\begin{equation} \tau^y_0 = \phi_2\big(\theta^y(z)\big) = -2\By_0\Gh_0 +\Bh_0\Gx_0, \end{equation}
\begin{equation} \tau^h_0 = \phi_2\big(\theta^h(z)\big) = -2\Bx_0\Gx_0 +2 \By_0\Gy_0. \end{equation}
Here $\phi_2$ denotes the projection $\Se_{(2)}\rightarrow \Se_{(2)}/\Se_{(1)}\subset gr(\Se)=P$. Define $\tau^u_k =  \partial^k \tau^u_0$, and let $P_{\tau}\subset P$ denote the subalgebra generated by the collection 
\begin{equation} \{\tau^u_k |\ k\ge 0,\ u=x,y,h\} .\end{equation} 
Since $\Theta\subset \Se^{\A_+}$, we have $\tau^u_0\in P^{\A_+}$. Since $P_{\tau}$ is generated as a $\partial$-ring by the $\tau^u_0$, it follows that $P_{\tau}\subset P^{\A_+}$.

The main result in this section is the following
\begin{thm} $P^{\A_+}=P_{\tau}$.
\end{thm}

{\it Proof of Theorem 1.3.}  Once Theorem 4.1 is established, Theorem 1.3 is an immediate consequence. Since
the vertex operators $\theta^u(z)$ already lie in $\Se^{\A_+}$, it
follows that the map $\Gamma: gr(\Se^{\A_+})\hookrightarrow
P^{\A_+}$ given by (3.9) is surjective, and hence is an isomorphism. By Lemma 3.6,
$\Se^{\A_+}=\Theta$. 

We will also show that
\begin{thm}  The map $\hat{\rho}: \mathcal{O}(sl(2),-K)\rightarrow \Se$, whose image is $\Theta$, is injective.\end{thm}
It follows that $\Se^{\A_+}$ is isomorphic to $\mathcal{O}(sl(2),-K)$, so we have a complete description of this commutant algebra.

\subsection{Outline of proof}
First, we will apply a classical theorem of Weyl to construct an isomorphism
$$\Phi:F/I\rightarrow P^{\A_0},$$ where $F$ is a polynomial
algebra on countably many variables corresponding to the quadratic generators of $P^{\A_0}$. By a linear change of variables, we may assume that the polynomials
$\tau^u_k$ correspond to a subset of the generators of $F$, which
generate a subalgebra $F_{T}\subset F$. We show that the
polynomials $\tau^u_k$ are algebraically independent, so that $F_T\bigcap I$ is trivial, and $F_T$ may be regarded as a subalgebra of $F/I$. Theorem 4.2 is an immediate consequence of this fact.

Using Lemma 3.12 (and assuming implicitly that we are working in
some $P_N$ for $N$ sufficiently large), we choose a monomial
ordering on $F$ and a corresponding Grobner basis for $I$ such that
elements of $F_{T}$ are in normal form. {\it By Lemma 3.13 and Remark 3.14, in order
to show that $P^{\A_+} = P_{\tau}$ it suffices to prove that
for any $\omega\in P^{\A_+}$ of positive degree, the normal form $\hat{\omega}\in
F$ of $\Phi^{-1}(\omega)$ contains a monomial in $F_{T}$ with nonzero coefficient.} 

$P$ has several $\mathbb{Z}_{\ge 0}$-gradings which will be useful to us. For a monomial
$$\mu = \Bx_{i_1}\cdots\Bx_{i_r}\By_{j_1}\cdots\By_{j_s}\Bh_{k_1}\cdots\Bh_{k_t} \Gx_{i'_1}\cdots\Gx_{i'_{r'}}\Gy_{j'_1}\cdots\Gy_{j'_{s'}}\Gh_{k'_1}\cdots\Gh_{k'_{t'}}\in P,$$
we define the $\be^u$-degree and $\ga^{u'}$-degree of $\mu$ as follows:
$$deg_{\Bx}(\mu) = r,\ \ \ \ \ deg_{\By}(\omega) = s,\ \ \ \ \ deg_{\Bh}(\mu) = t,$$ 
$$deg_{\Gx}(\mu) = r',\ \ \ \ \ deg_{\Gy}(\mu) = s',\ \ \ \ \ deg_{\Gh}(\mu) = t'.$$
Similarly, we define the $\be^u$-level and $\ga^{u'}$-level of $\mu$ to be
$$lev_{\Bx}(\mu) = \sum_{a=1}^r i_a,\ \ \ \ \ lev_{\By}(\omega) = \sum_{a=1}^s j_a,\ \ \ \ \ lev_{\Bh}(\mu) = \sum_{a=1}^t k_a,$$
$$lev_{\Gx}(\mu) = \sum_{a=1}^{r'} i'_a,\ \ \ \ \ lev_{\Gy}(\omega) = \sum_{a=1}^{s'} j'_a,\ \ \ \ \ lev_{\Gh}(\mu) = \sum_{a=1}^{t'} k'_a.$$

We will see that the condition $\omega\in P^{\A_+}$ implies that the projection of $\omega$ onto a certain homogeneous subspace (with respect to the above gradings) is nonzero. This will force $\hat{\omega}$ to contain a monomial in $F_T$ with non-zero coefficient.

\subsection{Description of $P^{\A_0}$}
It is immediate from (3.11)-(3.13) that as a module over
$\A_0=sl(2)$, $P$ is isomorphic to
\begin{equation} Sym\big(\bigoplus_{n\ge 0}W^1_n\oplus W^2_n\oplus W^3_n\big),\end{equation}
where each $W^i_n$ is a copy of the standard 2-dimensional irreducible
$sl(2)$-module. In particular, for each $n\ge 0$, each of the
following vector spaces form such a copy:
\begin{equation}W^1_n=\langle\Bx_n,\Gy_n\rangle,\ \ \
W^2_n=\langle\By_n,\Gx_n\rangle,\ \ \
W^3_n=\langle\Bh_n,\Gh_n\rangle. \end{equation} 

The description of the $sl(2)$-invariant subspace of such a module can be found in
\cite{W} (p. 45, Theorem 2.6.A, and p. 70, Theorem 2.14.A).
\begin{thm}
Let $S = Sym\big(\bigoplus_{n\ge 0} W_n\big)$, where $W_n = \langle
a^1_n,a^2_n\rangle$ is a copy of the standard 2-dimensional irreducible
$sl(2)$-module. The invariant subalgebra $S^{sl(2)}$ is generated by
the $2\times 2$-determinants:
\begin{equation}q_{ij}= \left|
\begin{array}{ll} a^1_i & a^2_i \\a^1_j & a^2_j
\end{array}\right| \ \ \ (0\le i<j),\end{equation}
which each correspond to a choice of two distinct modules from the
collection $\{W_n|\ n\ge 0\}$. The ideal of relations among the
polynomials $q_{ij}$ is generated by the polynomials
\begin{equation}r_{ijkl} = q_{ij}q_{kl}-q_{ik}q_{jl}+q_{il}q_{jk}.\end{equation} Each
of these polynomials corresponds to a choice of four distinct
modules from the collection $\{W_n|\ n\ge 0 \}$.
\end{thm}

In our context, taking into account the normalization of the modules
$\{W^i_n|\ i=1,2,3,\ n\ge 0\}$, $P^{\A_0}$ is generated by the
following six types of polynomials, which each correspond to a
choice of two distinct modules from the collection $\{W^i_n|\
i=1,2,3,\ n\ge 0\}$:
\begin{equation} q^{1,1}_{i,j} = 2\Bx_i\Gy_j -2\Bx_j\Gy_i,\ \ 0\le
i<j,\end{equation} 
\begin{equation} q^{2,2}_{i,j} = 2\By_i\Gx_j -2\By_j\Gx_i,\ \ 0\le
i<j,\end{equation} 
\begin{equation} q^{3,3}_{i,j} = \Bh_i\Gh_j -\Bh_j\Gh_i,\ \ 0\le
i<j,\end{equation} 
\begin{equation} q^{1,2}_{i,j} = 2\Bx_i\Gx_j
-2\By_j\Gy_i,\ \ i,j\ge 0,\end{equation} 
 \begin{equation} q^{1,3}_{i,j} = -2\Bx_i\Gh_j +
\Bh_j\Gy_i,\ \ i,j\ge 0,
\end{equation} 
\begin{equation} q^{2,3}_{i,j} = 2\By_i\Gh_j -
\Bh_j\Gx_i,\ \ i,j\ge 0. \end{equation}

Note that
\begin{equation}\tau^x_k = \sum_{i=0}^k
\binom{k}{i}q^{1,3}_{i,k-i},\ \ \ \tau^y_k = \sum_{i=0}^k
\binom{k}{i}q^{2,3}_{i,k-i},\ \ \ \tau^h_k = \sum_{i=0}^k
\binom{k}{i}q^{1,2}_{i,k-i}.\end{equation} It will be convenient to
perform a linear change of variables and replace $q^{1,3}_{0,k}$,
$q^{2,3}_{0,k}$, $q^{1,2}_{0,k}$ with $\tau^x_k$, $\tau^y_k$, and
$\tau^h_k$, respectively, using (4.16).

Let $F$ denote the polynomial algebra on the following generators: 
$$Q^{1,2}_{i,j}, Q^{1,3}_{i,j}, Q^{2,3}_{i,j},\ \ \  i>0,\ j\ge 0,$$
$$Q^{1,1}_{k,l},Q^{2,2}_{k,l}, Q^{3,3}_{k,l},\ \ \ 0\le k<l,$$
$$T^x_m, T^y_m, T^h_m,\ \ \ m\ge 0.$$ Let $I$ be the ideal
generated by the relations of the form (4.9). By Theorem 4.3, the
map \begin{equation} \Phi: F/I\rightarrow P^{\A_0}\end{equation} sending $Q^{a,b}_{i,j}\mapsto
q^{a,b}_{i,j}$ and $T^u_k\mapsto \tau^u_k$ is an isomorphism. Let
$F_T$ be the subalgebra of $F$ generated by the variables $T^u_k$.

\begin{rem} The description (4.6) of $P$ induces the auxiliary $\mathbb{Z}_{\ge 0}$-gradings $deg_{W^i}$ and $lev_{W^i}$ on $P$, defined as follows:
$$deg_{W^1} = deg_{\Bx}\ +\ deg_{\Gy},\ \ \ deg_{W^2} = deg_{\By}\ +\ deg_{\Gx},\ \ \ deg_{W^3} = deg_{\Bh}\ +\ deg_{\Gh},$$
$$lev_{W^1} = lev_{\Bx}\ +\ lev_{\Gy},\ \ \ lev_{W^2} = lev_{\By}\ +\ lev_{\Gx},\ \ \ lev_{W^3} = lev_{\Bh}\ +\ lev_{\Gh}.$$
$P^{\A_+}$ is graded by (total) level and $W^i$-degree, since each of the operators $v^u(n)$ for $u=x,y,h$ and $n\ge 0$ is homogeneous of level $-n$ and preserves $W^i$-degree. However, for $n>0$, $v^u(n)$ is not homogeneous with respect to $W^i$-level, so $P^{\A_+}$ is not graded by $W^i$-level. \end{rem}

\begin{rem} Since each $q^{a,b}_{i,j}$ and $\tau^u_k$ is homogeneous with respect to level and $W^i$-degree, $F$ and $I$ inherit these gradings in an obvious way. \end{rem}

Since the generators $q^{a,b}_{j,k}\in P^{\A_0}$ each correspond to a choice of two distinct modules from the collection $\{W^i_l|\ i=1,2,3,\ l\ge 0\}$ (namely $W^a_j$ and $W^b_k$), a monomial $$\mu = q^{a_1,b_1}_{j_1,k_1}\cdots q^{a_d,b_d}_{j_d,k_d}$$ of degree $d$ in the variables $q^{a,b}_{j,k}$ corresponds uniquely to the list of pairs \begin{equation} \mathcal{L}_{\mu} = \{\{W^{a_1}_{j_1},W^{b_1}_{k_1}\},\dots,\{W^{a_d}_{j_d},W^{b_d}_{k_d}\}\}.\end{equation} Consider the expansion of $\mu$ as a polynomial of degree $2d$ in the variables $\be^u_k,\ga^{u'}_k$. Each $q^{a,b}_{j,k}$ appearing in $\mu$ will contribute a factor of the form $\be\ga$, of which there are exactly two choices. (For example, each $q^{1,2}_{j,k}$ can contribute either $\Bx_j\Gx_k$ or $\By_k\Gy_j$). 

Suppose for the moment that we consider only monomials $\mu$ in the variables $q^{a,b}_{j,k}$ for $a<b$, which are given by (4.13)-(4.15). 

\begin{lem} Let $$\epsilon = (\Bx_0\Gh_{i_1})\cdots (\Bx_0\Gh_{i_r} )(\Bx_0\Gx_{j_1})\cdots (\Bx_0\Gx_{j_s})( \By_{k_1}\Gh_0)\cdots (\By_{k_t}\Gh_0).$$ The only monomials in the variables $q^{a,b}_{j,k}$ for $a<b$ which contain $\epsilon$ with nonzero coefficient are of the form  $$\mu=q^{1,3}_{0,i'_1}\cdots q^{1,3}_{0,i'_r} q^{1,2}_{0,j_1}\cdots q^{1,2}_{0,j_s} q^{2,3}_{k_1,i''_1}\cdots q^{2,3}_{k_t,i''_t},$$ where the lists $(i'_1,\dots,i'_r,i''_1,\dots,i''_t)$ and $(i_1,\dots,i_r,0,\dots,0)$ are related by a permutation.
\end{lem}
\begin{proof} First, the only variables $q^{a,b}_{j,k}$ for $a<b$ which can contribute $\Bx_0$ are $q^{1,2}_{0,k}$ (which contains the monomial $\Bx_0\Gh_k$) and $q^{1,3}_{0,k}$ (which contains $\Bx_0\Gx_0$). Since $\epsilon$ is divisible by $(\Bx_0)^{r+s}$, exactly $r+s$ of the variables $q^{1,2}_{0,k},q^{1,3}_{0,k}$ must occur. Since $\mu$ is not divisible by any $q^{2,2}_{j,k}$, no pairings of the form $\By_{k_a}\Gx_{j_b}$ can occur, so each $\Gx_{j_b}$ must be paired with one of the $\Bx_0$'s, for $b=1,\dots,s$. Hence $\mu$ must be divisible by exactly $s$ of the variables $q^{1,2}_{0,k}$, namely $q^{1,2}_{0,j_1}\cdots q^{1,2}_{0,j_s}$. Then $\mu$ must be divisible by exactly $r$ of the variables $q^{1,3}_{0,k}$, say  $q^{1,3}_{0,i'_1}\cdots q^{1,3}_{0,i'_r}$. The indices $i'_1,\dots,i'_r$ correspond to a choice of $r$ modules $\{W^3_{i'_1},\dots,W^3_{i'_r}\}$, from the set $\{W^3_{i_1},\dots,W^3_{i_r},W^3_0,\dots W^3_0\}$, which contain $r+t$ elements.

Let $\{W^3_{i''_1},\dots,W^3_{i''_t}\}\subset\{W^3_{i_1},\dots,W^3_{i_r},W^3_0,\dots W^3_0\}$ be the complement of $\{W^3_{i'_1},\dots,W^3_{i'_r}\}$, so that $(i'_1,\dots,i'_r,i''_1,\dots,i''_t)$ is some permutation of $(i_1,\dots,i_r,0,\dots,0)$. The factor $$q^{1,3}_{0,i'_1}\cdots q^{1,3}_{0,i'_r}q^{1,2}_{0,j_1}\cdots q^{1,2}_{0,j_s}$$ appearing in $\mu$ accounts for the factor $(\Bx_0\Gh_{i'_1})\cdots (\Bx_0\Gh_{i'_r} )(\Bx_0\Gx_{j_1})\cdots (\Bx_0\Gx_{j_s})$ appearing in $\epsilon$.  The remaining factor $(\By_{k_1}\Gh_{i''_1})\cdots (\By_{k_t}\Gh_{i''_t})$ of $\epsilon$ can only appear in a monomial in the variables $q^{a,b}_{j,k}$ of the form $q^{2,3}_{k_1,i''_1}\cdots q^{2,3}_{k_t,i''_t}$.
\end{proof}

\begin{lem} As above, let $$\epsilon = (\Bx_0\Gh_{i_1})\cdots (\Bx_0\Gh_{i_r} )(\Bx_0\Gx_{j_1})\cdots (\Bx_0\Gx_{j_s})( \By_{k_1}\Gh_0)\cdots (\By_{k_t}\Gh_0).$$ 
The only monomials in the variables $\tau^u_k$ which can contain $\epsilon$ with nonzero coefficient are of the form
$$\nu=\tau^x_{i'_1}\cdots\tau^x_{i'_r}\tau^h_{j_1}\cdots\tau^h_{j_s}\tau^y_{k'_1}\cdots\tau^y_{k'_t},$$ where the lists $(i'_1,\dots,i'_r)$ and $(k'_1,\dots,k'_t)$ are obtained from the lists $(i_1,\dots,i_r)$ and $(k_1,\dots,k_t)$ by replacing some of the pairs $(i_a,k_b)$ with $(0,i_a+k_b)$. 
\end{lem}

\begin{proof} In order for $\nu$ to contain $\epsilon$ with nonzero coefficient, a monomial of the form
$$\mu=q^{1,3}_{0,i'_1}\cdots q^{1,3}_{0,i'_r} q^{1,2}_{0,j_1}\cdots q^{1,2}_{0,j_s} q^{2,3}_{k_1,i''_1}\cdots q^{2,3}_{k_t,i''_t},$$
must appear when $\nu$ is expanded as a polynomial in the variables $q^{a,b}_{j,k}$, using (4.16). Here $(i'_1,\dots,i'_r,i''_1,\dots,i''_t)$ is some permutation of $(i_1,\dots,i_r,0,\dots,0)$, by Lemma 4.6. It is immediate from (4.16) that the only monomial in the variables $\tau^u_k$ which will contain $\mu$ is
$$\tau^x_{i'_1}\cdots\tau^x_{i'_r}\tau^h_{j_1}\cdots\tau^h_{j_s}\tau^y_{i''_1+k_1}\cdots\tau^y_{i''_t+k_t}.$$ Setting $k'_a = i''_a+k_a$ for each $a=1,\dots,t$, the claim follows.
\end{proof}

\begin{lem} The polynomials $\tau^u_k$ are algebraically
independent. Equivalently, $F_T\bigcap I$ is trivial, so we may regard $F_T$ as a subalgebra of $F/I$, which maps isomorphically onto $P_{\tau}$ under $\Phi$.
\end{lem}

\begin{proof}
Let $Q$ be a polynomial in the variables $\tau^u_k$, which we may
assume to be homogeneous of fixed level and $W^i$-degree. Let
$$\mu=\tau^x_{i_1}\cdots\tau^x_{i_r}\tau^h_{j_1}\cdots\tau^h_{j_s}\tau^y_{k_1}\cdots\tau^y_{k_t}$$
be a monomial appearing in $Q$. Clearly
$$deg_{W^1}(\mu)=r+s,\ \ \ deg_{W^2}(\mu)=s+t,\ \ \ deg_{W^3}(\mu)=r+t.$$
Since $Q$ is homogeneous with respect to $W^1$-degree, $W^2$-degree,
and $W^3$-degree, it follows that $deg_{W^i}(\mu)=deg_{W^i}(Q)$ for
$i=1,2,3$. Solving for $r$, $s$, and $t$, we obtain:

$$r=\frac{1}{2}\big(deg_{W^1}(Q)-deg_{W^2}(Q)+deg_{W^3}(Q)\big),$$
$$s=\frac{1}{2}\big(deg_{W^1}(Q)+deg_{W^2}(Q)-deg_{W^3}(Q)\big),$$
$$t=\frac{1}{2}\big(-deg_{W^1}(Q)+deg_{W^2}(Q)+deg_{W^3}(Q)\big).$$
Since $r$, $s$, and $t$ only depend on $deg_{W^i}(Q)$, for
$i=1,2,3$, they are the same for all monomials $\mu$ appearing in
$Q$.

Fix a monomial
$\mu=\tau^x_{i_1}\cdots\tau^x_{i_r}\tau^h_{j_1}\cdots\tau^h_{j_s}\tau^y_{k_1}\cdots\tau^y_{k_t}$
appearing $Q$ such that the number of zeros appearing in the list
$\{i_1,\dots,i_r\}$ is maximal (in the case $r=0$, no such choice is necessary). By Lemma 4.7, the monomial $$\epsilon=(\be^x_0\Gh_{i_1})\cdots (\be^x_0\Gh_{i_r})(\be^x_0\Gx_{j_1})
\cdots (\Bx_0\Gx_{j_s})(\By_{k_1}\Gh_0)\cdots (\By_{k_t}\Gh_0)$$
appears in $\mu$ with nonzero coefficient.  Moreover, any other monomial containing $\epsilon$ with nonzero coefficient has the form
$$\mu'=\tau^x_{i'_1}\cdots\tau^x_{i'_r}\tau^h_{j_1}\cdots\tau^h_{j_s}\tau^y_{k'_1}\cdots\tau^y_{k'_t},$$
where the lists $(i'_1,\dots,i'_r)$ and $(k'_1,\dots,k'_t)$ are obtained from $(i_1,\dots,i_r)$ and $(k_1,\dots,k_t)$ by replacing some of the pairs $(i_a,k_b)$ with $(0,i_a+k_b)$. Since the number of zeros in the list $\{i_1,\dots,i_r\}$ is maximal, no such $\mu'$ can appear in $Q$.
Hence $\epsilon$ appears in $Q$, and in particular, $Q\neq 0$.
\end{proof}

\begin{cor} $\hat{\rho}:\mathcal{O}\big(sl(2),-K\big)\rightarrow
\Se$ is injective.
\end{cor}
\begin{proof} Recall that $\hat{\rho}$ is a morphism in the category $\Re$ if
we declare that the generators $x(z),y(z),h(z)$ of
$\mathcal{O}(sl(2),-K)$ have degree 2. The
induced map on the associated graded algebras
$$gr(\hat{\rho}):gr\big(\mathcal{O}(sl(2),-K)\big)\rightarrow P$$
is a $\partial$-ring homomorphism. By Corollary 3.4,
$gr\big(\mathcal{O}(sl(2),-K)\big)$ is the
polynomial algebra with generators $x_k,y_k,h_k$ for $k\ge 0$, and $gr(\hat{\rho})$ sends $$x_k\mapsto
\tau^x_k,\ \ \ \ y_k\mapsto \tau^y_k,\ \
\ \ h_k\mapsto\tau^h_k.$$ Since the polynomials $\tau^u_k$ are
algebraically independent, it follows that $gr(\hat{\rho})$ is
injective, so $\hat{\rho}$ must be injective as well. \end{proof}

Next, we will choose a monomial ordering on $F$ and a corresponding Grobner basis $B$ for $I$. We order the generators
$Q^{a,b}_{i,j},T^u_k$ as follows:
\begin{equation} Q^{3,3}_{i,j}> Q^{2,3}_{i,j}> Q^{2,2}_{i,j}>Q^{1,2}_{i,j}>Q^{1,1}_{i,j}>T^y_k>T^x_k>T^h_k,\end{equation}
for all $i,j,k$. Then for each $a,b = 1,2,3$ such that
$a\le b$,
\begin{equation}Q^{a,b}_{i,j}> Q^{a,b}_{k,l}\end{equation} if $j>l$ or if $j=l$ and $i>k$. Likewise, $T^u_k>T^u_l$ if $k>l$, for each $u=x,y,h$.
Finally, order monomials in these variables using the standard
lexicographic ordering. To find a Grobner basis for $I$, begin with
the generating set of relations of the form (4.9), eliminating the variables
$q^{1,2}_{0,k}$, $q^{1,3}_{0,k}$, and $q^{2,3}_{0,k}$ using (4.16). Replacing $q^{a,b}_{i,j},\tau^u_k$ with $Q^{a,b}_{i,j},T^u_k$, respectively, we obtain a generating set $B'$ for $I$.
Extend this set to a Grobner basis $B$ for $I$ using Buchberger's
algorithm. For any $\omega\in P^{\A_0}$, let $\hat{\omega}\in F$
denote the corresponding normal form. 

\begin{rem} It follows from Lemmas 3.12 and 4.8, and the term ordering (4.19)-(4.20) that any element of $F_T$ is automatically in normal form with respect to $B$. \end{rem}

The following observation will be useful to us later:

\begin{lem} Suppose that $f\in F$ is in normal form. Then $f$ is not
divisible by the product $Q^{3,3}_{i,j}T^h_0$ for any $0\le
i<j$. Similarly, $f$ is not divisible by $Q^{2,3}_{i,j}
T^h_0$ for any $i>0$ and $j\ge 0$.
\end{lem}
\begin{proof}
We begin with the first statement. It suffices to show that
$Q^{3,3}_{i,j}T^h_0$ is the leading term of an element of $B$. Note
that \begin{equation} q^{3,3}_{i,j} q^{1,2}_{0,0}-q^{1,3}_{0,i}q^{2,3}_{0,j} +
q^{1,3}_{0,j}q^{2,3}_{0,i}\end{equation} is a relation of the form (4.10). Use (4.16) to eliminate the variables $q^{1,2}_{0,0}$, $q^{1,3}_{0,i}$, $q^{1,3}_{0,j}$, $q^{2,3}_{0,i}$, and $q^{2,3}_{0,j}$ from (4.21), and then replace $q^{a,b}_{r,s}$, $\tau^u_t$ with $Q^{a,b}_{r,s}$, $T^u_{t}$, respectively. We obtain the relation
$$Q^{3,3}_{i,j} T^h_0 - \bigg(T^x_i -
\sum_{a=1}^i \binom{i}{a}Q^{1,3}_{a,i-a}\bigg)\bigg(T^y_j - \sum_{b=1}^j \binom{j}{b}Q^{2,3}_{b,j-b}\bigg) + $$ \begin{equation} \bigg(T^x_j - \sum_{a=1}^j
\binom{j}{a}Q^{1,3}_{a,j-a}\bigg) \bigg(T^y_i -
\sum_{b=1}^i \binom{i}{b}Q^{2,3}_{b,i-b}\bigg),\end{equation} which lies in $B'$ (and hence in $B$)
by definition. According to the term ordering given by (4.19)-(4.20), the leading term of (4.22) is $Q^{3,3}_{i,j}T^h_0$, as desired. The proof of the second statement in Lemma 4.11 is verbatim.
\end{proof}

\subsection{Some commutative algebra}
Since $P^{\A_0}$ and $P^{\A_+}$ are subalgebras of the polynomial algebra $P$, they are integral domains. Hence $F/I$ is also a domain, as is any subalgebra of $F/I$. For any domain $R$, we shall denote the corresponding field of fractions by $\overline{R}$. The isomorphism $\Phi:F/I\rightarrow P^{\A_0}$ given by (4.17) extends to an isomorphism
\begin{equation}\overline{\Phi}: \overline{F/I} \rightarrow \overline{P^{\A_0}}\subset \overline{P}.\end{equation} If $S$ is any multiplicative subset of $F/I$, we may regard the localization $S^{-1}(F/I)$ as a subalgebra of $\overline{F/I}$, and by (4.23), as a subalgebra of $\overline{P}$. Recall that $F_T$ may be regarded as a subalgebra of $F/I$, and that $\Phi$ maps $F_T$ isomorphically onto $P_{\tau}$, by Lemma 4.8. 

For the remainder of this section, let $S\subset F_T$ be the multiplicative subset generated by the single element $T^h_0$. Let $R$ be the subalgebra $\Phi^{-1}(P^{\A_+})\subset F/I$. Clearly $F_T\subset R$; our goal is to prove that $F_T= R$. Since $S\subset F_T\subset R\subset F/I$, we may localize all these rings with respect to $S$. We obtain inclusions 
$$ S^{-1} (F_T)\hookrightarrow S^{-1} R\hookrightarrow S^{-1} (F/I)\hookrightarrow\Phi(S)^{-1} P.$$ Here $\Phi(S)$ denotes the multiplicative subset of $P$ generated by $\Phi(T^h_0) = \tau^h_0$, and the last inclusion above is the restriction of $\overline{\Phi}$ to the subalgebra $ S^{-1} (F/I)\subset \overline{F/I}$. We need the following technical statement.

\begin{lem}
Let $\alpha\in S^{-1} (F_T)$. If $\overline{\Phi}(\alpha)\in P$, then $\alpha\in F_T$. \end{lem}

\begin{proof}
Equivalently, we need to prove that $\alpha\notin F_T$ implies that $\overline{\Phi}(\alpha)\notin P$. Let us rephrase this as an ideal membership problem. Since $F_T$ is a polynomial algebra, $\alpha$ can be written uniquely in the form
$$\alpha =  \sum_{i\ge 0} \frac{\alpha_i}{(T^h_0)^i},$$ where the $\alpha_i$'s are elements of $F_T$ which are not divisible by $T^h_0$ for $i>0$ (and hence do not lie in the principal ideal $J\subset F_T$ generated by $T^h_0$). The condition $\alpha\notin F_T$ means that some $\alpha_i\neq 0$ for $i>0$.
By multiplying by an appropriate power of $T^h_0$, we may assume without loss of generality that 
$$\alpha =\alpha_0 +\frac{\alpha_1}{T^h_0},$$ where $\alpha_1\notin J$.

The statement $\overline{\Phi}(\alpha)\notin P$ is equivalent to the statement that $\Phi(\alpha_1)$ does not lie in the ideal $\mathcal{J}\subset P$ generated by $\tau^h_0 = \Phi(T^h_0)$. 
Without loss of generality, we may assume that $\alpha_1$ is a linear combination of monomials of the form
$$\mu = T^x_{i_1}\cdots T^x_{i_r} T^h_{j_1}\cdots T^h_{j_s} T^y_{k_1}\cdots T^y_{k_t}.$$
for which $r,s,t$ are fixed. Moreover, we may assume that each such monomial $\mu\notin J$, so that $(j_1,\dots,j_s)$ are all positive.

Recall from the proof of Lemma 4.8 that for each monomial $$\mu = T^x_{i_1}\cdots T^x_{i_r} T^h_{j_1}\cdots T^h_{j_s} T^y_{k_1}\cdots T^y_{k_t}$$ appearing in $\alpha_1$ for which the list $(i_1,\dots,i_r)$ contains the maximum number of zeros, $\Phi(\alpha_1)$ will contain the monomial
$$\epsilon = (\Bx_0\Gh_{i_1})\cdots (\Bx_0\Gh_{i_r})( \Bx_0\Gx_{j_1})\cdots (\Bx_0\Gx_{j_s})( \By_{k_1}\Gh_0)\cdots (\By_{k_t}\Gh_0)$$ with nonzero coefficient.

Since $\tau^h_0 =-2\Bx_0\Gx_0+2\By_0\Gy_0$, any element $\omega\in \mathcal{J}$ has the property that each monomial appearing in $\omega$ is divisible by either $\Bx_0\Gx_0$ or $\By_0\Gy_0$. Since $\epsilon$ is not divisible by either $\Gx_0$ (since $j_1,\dots,j_s$ are all positive), or $\Gy_0$, we conclude that $\Phi(\alpha_1)\notin \mathcal{J}$, as claimed.
\end{proof}

\begin{cor} Let $\omega\in P^{\A_+}$. If $\tau^h_0\omega\in P_{\tau}$, then $\omega\in P_{\tau}$ as well. \end{cor}

\begin{proof} 
The condition $\tau^h_0\omega\in P_{\tau}$ means that $\Phi^{-1}(\tau^h_0\omega)$ can be expressed as a polynomial $\nu\in F_T$. Consider the element $\frac{1}{T^h_0}\nu\in S^{-1}(F_T)$ and note that
$$\overline{\Phi}\big(\frac{1}{T^h_0}\nu\big) = \omega \in P.$$ It follows from Lemma 4.12 that $\frac{1}{T^h_0}\nu\in F_T$, so that $\nu = T^h_0 \nu'$ for some $\nu'\in F_T$. Then
$$\omega = \overline{\Phi}\big(\frac{1}{T^h_0} \nu\big) = \Phi(\nu') $$ which lies in $P_{\tau}$ since $\nu'\in F_T$. It follows that $\omega\in P_{\tau}$, as claimed. \end{proof}

By applying Corollary 4.13 repeatedly, we see that for any $\omega\in P^{\A_+}$ and $r\ge 0$, \begin{equation}(\tau^h_0)^r\omega \in P_{\tau}  \Longleftrightarrow \omega\in P_{\tau}.\end{equation}Thus given $\omega\in P^{\A_+}$, in order to prove that $\omega\in P_{\tau}$, it suffices to prove that $(\tau^h_0)^r\omega\in P_{\tau}$ for some $r$.

\begin{lem} If $\omega\in P^{\A_+}$ is homogeneous of degree $2d$, choose $r>2d$, and let $\omega' = (\tau^h_0)^r\omega$. Then each monomial appearing in the normal form $\widehat{\omega'}$ is divisible by $T^h_0$. \end{lem}
\begin{proof} Recall that each generator $q^{a,b}_{j,k}\in P^{\A_0}$ corresponds to the pair $\{W^a_j,W^b_k\}$ of distinct modules from the collection $\{W^i_l|\ i=1,2,3,\ l\ge 0\}$. By Theorem 4.3, we may write $\omega$ as a polynomial $\tilde{\omega}$ of degree $d$ in the variables $q^{a,b}_{j,k}$. Recall that each monomial 
$$\mu = q^{a_1,b_1}_{j_1,k_1}\cdots q^{a_d,b_d}_{j_d,k_d}$$ appearing in $\tilde{\omega}$, corresponds to the list of pairs modules $$\mathcal{L}_{\mu} = \{\{W^{a_1}_{j_1},W^{b_1}_{k_1}\},\dots,\{W^{a_d}_{j_d},W^{b_d}_{k_d}\}\}.$$ If we choose $r>2d$, the corresponding list 
$$ \mathcal{L}_{(\tau^h_0)^r\mu}=\{\{W^1_0,W^2_0\},\dots,\{W^1_0, W^2_0\},\{W^{a_1}_{j_1},W^{b_1}_{k_1}\},\dots,\{W^{a_d}_{j_d},W^{b_d}_{k_d}\}\}$$  will contain at least $2d+1$ copies of $W^1_0$ and $2d+1$ copies of $W^2_0$. Any monomial $\mu'$ in the variables $q^{a,b}_{j,k}$ for which $\mathcal{L}_{\mu'}$ and $\mathcal{L}_{(\tau^h_0)^r\mu}$ contain the same collection of modules (possibly re-ordered) must be divisible by $\tau^h_0$ by the pigeonhole principle, since each $q^{a,b}_{j,k}$ depends on a pair of {\it distinct} modules. This statement remains true after making the change of variables (4.16), so in particular the normal form $\hat{\omega}$ will have the desired property.
\end{proof}

\subsection{Description of $P^{\A_+}$}
Let $\omega\in P^{\A_+}$ be nonzero, which we may assume to be homogeneous of
level $l$ and $W^i$-degree $d_i$ (and hence total degree $ d = d_1+d_2+d_3$). Note that $\tau^h_0$ is homogeneous of level $0$ and $W^1$-degree, $W^2$-degree, $W^3$-degree $1,1,0$, respectively, so $(\tau^h_0)^r\omega$ will still be homogeneous with respect to these gradings. By (4.24) and Lemma 4.14, we may assume without loss of generality that each monomial appearing in the normal form $\hat{\omega}$ is divisible by $T^h_0$. In particular, $d_1>0$ and $d_2>0$. We will show that the condition $\omega\in P^{\A_+}$ implies that the projection of $\omega$ onto a certain homogeneous subspace is nonzero. This will force $\hat{\omega}$ to contain a monomial in $F_T$ with nonzero coefficient. Theorem 4.1 then follows immediately from Lemma 3.13 and Remark 3.14.

\begin{lem} Let $\pi^{\Gy,e}$ denote the projection of $P$ onto its homogeneous component of $\Gy$-degree $e$. Then $\pi^{\Gy,0}(\omega)\neq 0$. Equivalently, $\omega$ has a nonzero term of $\Bx$-degree $d_1$, since $deg_{W^1}=deg_{\Bx}+ deg_{\Gy}$ and $deg_{W^1}(\omega)=d_1$. 
\end{lem}

\begin{proof} Let $e$ be the minimal $\Gy$-degree of terms appearing
in $\omega$, and write $\omega = \omega_0 +\omega_1$, where $\omega_0 = \pi^{\Gy,e}(\omega)$.
Let $k$ be the largest
integer such that $\Gy_k$ appears in $\omega_0$, and write
$$\omega_0 =  (\Gy_k)^t q_t + (\Gy_k)^{t-1}q_{t-1}+\dots+
\Gy_k q_1 + q_0,$$ where the $q_a$'s do not depend on $\Gy_k$, and at
least one of the $q_a$'s is nonzero for $a=1,\dots,t$. Clearly each nonzero
$q_a$ must have $\Gy$-degree $e-a$ since $\omega_0$ has $\Gy$-degree $e$. 

Applying (3.13) in the case $\G=sl(2)=V$, and working in the basis (4.1), we have 
\begin{equation} v^y(n)(\Gy_m) = -\frac{1}{2} m(m-1)\cdots (m-n+1)\Bx_{m-n},\end{equation}
and in particular, $v^y(k)(\Gy_k) = -\frac{1}{2}k! \Bx_0$. We compute
$$v^y(k)(\omega) = t (\Gy_k)^{t-1}(-\frac{1}{2}k! \Bx_0)q_t+ (t-1) (\Gy_k)^{t-2}(-\frac{1}{2}k! \Bx_0)q_{t-1}+\dots + (-\frac{1}{2}k! \Bx_0)q_1$$  $$+ (\Gy_k)^t v^y(k)(q_t) +
(\Gy_k)^{t-1} v^y(k)(q_{t-1})+\dots + (\Gy_k) v^y(k) (q_1) +
v^y(k)(q_0) $$ $$+ v^y(k)(\omega_1).$$ The term
$$t (\Gy_k)^{t-1}(-\frac{1}{2}k! \Bx_0)q_t + (t-1)  (\Gy_k)^{t-2}(-\frac{1}{2}k!
\Bx_0)q_{t-1} +\dots + (-\frac{1}{2}k! \Bx_0)q_1$$ is homogeneous of $\Gy$-degree $e-1$. Furthermore, this term is nonzero since at
least one of the $q_a$'s is nonzero, and none of the $q_a$'s depends
on $\Gy_k$. 

We claim that
\begin{equation}\pi^{\Gy,e-1}\big(v^y(k)(\omega)\big)
=\end{equation} $$t (\Gy_k)^{t-1}(-\frac{1}{2}k! \Bx_0)q_t + (t-1) (\Gy_k)^{t-2}(-\frac{1}{2}k! \Bx_0)q_{t-1} +\dots + (-\frac{1}{2}k!
\Bx_0)q_1.$$ In other words, none of
the other terms appearing in $v^y(k)(\omega)$ can have
$\Gy$-degree $e-1$. It follows that $v^y(k)(\omega)\neq 0$, which contradicts $\omega\in P^{\A_+}$. Hence we must have $e=0$.

In order to prove (4.26), we need to show that 
$$(\Gy_k)^t v^y(k)(q_t) +
(\Gy_k)^{t-1} v^y(k)(q_{t-1})+\dots + (\Gy_k) v^y(k) (q_1) +
v^y(k)(q_0) + v^y(k)(\omega_1)$$
has no term of $\Gy$-degree $e-1$. By (4.25), $v^y(k)(\Gy_m)= 0$ for $m<k$, and $q_a$ does not depend on $\Gy_m$ for any $m\ge
k$. It follows that each term of $v^y(k)(q_a)$ must have $\Gy$-degree $e-a$, so
each term appearing in
$$(\Gy_k)^t v^y(k)(q_t) +
(\Gy_k)^{t-1} v^y(k)(q_{t-1})+\dots + (\Gy_k) v^y(k) (q_1)+
v^y(k)(q_0)$$ has $\Gy$-degree $e$.

Finally, $\omega_1$ consists of terms with $\Gy$-degree at least
$e+1$, so every term of $v^y(k)(\omega_1)$ has $\Gy$-degree at least
$e$. \end{proof}

\begin{lem} Let $\pi^{\Bx}_e$ denote the projection of $P$ onto its homogeneous component of $\Bx$-level $e$. Then $\pi^{\Bx}_0\circ \pi^{\Gy,0}(\omega)\neq 0$. Equivalently, $\omega$ has a nonzero term which does not depend on any of the variables $\Gy_i$ for $i\ge 0$ and $\Bx_j$ for $j>0$.
\end{lem}

\begin{proof} Write $\omega = \omega_0+\omega_1$, where $\omega_0 = \pi^{\Gy,0}(\omega)$. Recall that $\omega_0$ has $\Bx$-degree $d_1>0$. Let $e$ be the minimal $\Bx$-level of terms appearing in $\omega_0$, and write $\omega_0 = \omega_0'+\omega_0''$, where $$\omega_0' = \pi^{\Bx}_e(\omega_0) = \pi^{\Bx}_e\circ \pi^{\Gy,0}(\omega).$$

Suppose that $e>0$, and let $k$ be the largest integer such that $\Bx_k$ appears in $\omega_0'$. Clearly $0<k\le e$. Write
$$\omega_0' =  (\Bx_k)^t q_t + (\Bx_k)^{t-1}q_{t-1}+\dots+
\Bx_k q_1 + q_0,$$ where the $q_a$'s do not depend on $\Bx_k$, and at
least one of the $q_a$'s is nonzero for $a=1,\dots,t$. Clearly each nonzero
$q_a$ must have $\Bx$-level $e-ka$, since $\omega_0'$ has $\Bx$-level $e$. 

By (3.11), we have $$v^h(n)\Bx_m = -m(m-1)\cdots (m-n+1)\Bx_0,$$ and in particular, $v^h(k)(\Bx_k) = -k! \Bx_0$. We compute
$$v^h(k)(\omega) = t (\Bx_k)^{t-1}(-k! \Bx_0)q_t+ (t-1) (\Bx_k)^{t-2}(-k! \Bx_0)q_{t-1}+\dots + (-k! \Bx_0)q_1$$  $$+ (\Bx_k)^t v^h(k)(q_t) +
(\Bx_k)^{t-1} v^h(k)(q_{t-1})+\dots + (\Bx_k) v^h(k) (q_1) +
v^h(k)(q_0) $$ $$+ v^h(k)(\omega_0'') + v^h(k)(\omega_1).$$ 

We claim that $$\pi^{\Bx}_{e-k}\circ\pi^{\Gy,0}\big(v^h(k)(\omega)\big) = $$ $$ t (\Bx_k)^{t-1}(-k! \Bx_0)q_t + (t-1)  (\Bx_k)^{t-2}(-k!
\Bx_0)q_{t-1} +\dots + (-k! \Bx_0)q_1, $$ which is clearly nonzero, contradicting $\omega\in P^{\A_+}$. This proves that $e=0$.

First, the term 
$$(\Bx_k)^t v^h(k)(q_t) +
(\Bx_k)^{t-1} v^h(k)(q_{t-1})+\dots + (\Bx_k) v^h(k) (q_1) +
v^h(k)(q_0) $$ must be homogeneous of $\Bx$-level $e$, since $v^h(k)(\Bx_m) = 0$ for $m<k$, and $q_a$ does not depend on $\Bx_m$ for $m\ge k$.

Second, the term $v^h(k)(\omega_0'')$ must have $\Bx$-level at least $e+1-k$ since $\omega_0''$ has $\Bx$-level at least $e+1$ and $v^h(k)$ can lower the $\Bx$-level by at most $k$.

Finally, the term $v^h(k)(\omega_1)$ must have positive $\Gy$-degree, since $\omega_1$ has positive $\Gy$-degree and $v^h(k)$ preserves $\Gy$-degree. \end{proof}

\begin{lem} Let $\pi^{\By,e}$ denote the projection of $P$ onto its homogeneous component of $\By$-degree $e$. Then $\pi^{\By,0}\circ\pi^{\Bx}_0\circ\pi^{\Gy,0}(\omega)\neq 0$. \end{lem}

\begin{proof} Let $\omega_0 = \pi^{\Bx}_0\circ\pi^{\Gy,0}(\omega)$, and write $\omega = \omega_0+\omega_1$. Let $e$ be the minimal $\By$-degree of terms appearing
in $\omega_0$, write $\omega_0 = \omega_0'+\omega_0''$, where
$$\omega_0' = \pi^{\By,e}(\omega_0) = \pi^{\By,e}\circ \pi^{\Bx}_0\circ\pi^{\Gy,0}(\omega).$$

Suppose that $e>0$, and let $k$ be the maximum integer such that $\By_k$ appears in
$\omega_0'$. Write
$$\omega_0' =  (\By_k)^t q_t + (\By_k)^{t-1}q_{t-1}+\dots+
\By_k q_1 + q_0,$$ where the $q_a$'s do not depend on $\By_k$, and at
least one of the $q_a$'s is nonzero for $a=1,\dots,t$. Clearly each nonzero
$q_a$ must have $\By$-degree $e-a$ since $\omega_0'$ has $\By$-degree $e$. 

By (3.12), we have $v^x(n)(\By_m) = -\frac{1}{2}m(m-1)\cdots(m-n+1)\Gx_{m-n}$, and in particular, $v^x(k)(\By_k) = -\frac{1}{2} k! \Gx_0$. We compute
$$v^x(k)(\omega) =  t (\By_k)^{t-1}(-\frac{1}{2}k! \Gx_0)q_t+ (t-1) (\By_k)^{t-2}(-\frac{1}{2}k! \Gx_0)q_{t-1}+\dots + (-\frac{1}{2}k! \Gx_0)q_1$$  $$+ (\By_k)^t v^x(k)(q_t) +
(\By_k)^{t-1} v^x(k)(q_{t-1})+\dots + (\By_k) v^x(k) (q_1) +
v^x(k)(q_0) $$ $$+ v^x(k)(\omega_0'') + v^x(k)(\omega_1).$$ 

We claim that 
$$\pi^{\By,e-1}\circ\pi^{\Bx}_0\circ\pi^{\Gy,0}\big(v^h(k)\omega)\big) = $$ $$ t (\By_k)^{t-1}(-\frac{1}{2}k! \Gx_0)q_t+ (t-1) (\By_k)^{t-2}(-\frac{1}{2}k! \Gx_0)q_{t-1}+\dots + (-\frac{1}{2}k! \Gx_0)q_1,$$ which is clearly nonzero, contradicting $\omega\in P^{\A_+}$. This proves that $e=0$.

First, the term $$(\By_k)^t v^x(k)(q_t) +
(\By_k)^{t-1} v^x(k)(q_{t-1})+\dots + (\By_k) v^x(k) (q_1) +
v^x(k)(q_0)$$ must have $\By$-degree $e$, since $q_t$ does not depend on $\By_m$ for any $m\ge k$.

Second, the term $v^x(k)(\omega_0'')$ must have $\By$ degree at least $e$, since $\omega_0''$ has $\By$-degree at least $e+1$.

Finally, the term $v^x(k)(\omega_1)$ must have positive $\Bx$-level or positive $\Gy$-degree. This follows from the fact that $$v^x(k)(\Bx_m) = -\frac{1}{2} m(m-1)\cdots(m-k+1)\Gy_{m-k,}$$ which shows that $v^x(k)$ can only lower the $\Bx$-level by raising the $\Gy$-degree.
\end{proof}

{\it Proof of Theorem 4.1.}
In order to prove Theorem 4.1, it suffices to show that the normal form $\hat{\omega}$ contains a monomial in $F_T$ with nonzero coefficient, by Lemma 3.13 and Remark 3.14. As we shall see, Lemmas 4.17 and 4.11 taken together will force $\hat{\omega}$ to contain such a monomial.

Let $\omega_0=\pi^{\By,0}\circ\pi^{\Bx}_0\circ \pi^{\Gy,0}(\omega)$, which is non-zero by Lemma 4.17. Since $\omega_0$ only depends on the variables $\Bx_0,\Gx_k,\Bh_k,\Gh_k$ for $k\ge 0$, it is a linear combination of monomials of the form
$$\epsilon=(\Bx_0)^{d_1}\Gx_{i_1}\cdots \Gx_{i_r}\Bh_{j_1}\cdots\Bh_{j_s}\Gh_{k_1}\cdots\Gh_{k_t}.$$
Here $r=d_2$, $s+t = d_3$, and $d_1+s = d_2+t = \frac{d}{2}$, since the total number of $\be$ and $\gamma$ are equal. A similar argument to the proof of Lemma 4.6 shows that the only possible monomials in the variables $q^{a,b}_{i,j}$ (before making the change of variables given by (4.16)) which can contain $\epsilon$ are of the form
$$q^{1,3}_{0,k'_1}\cdots q^{1,3}_{0,k'_a} q^{1,2}_{0,i'_1}\cdots q^{1,2}_{0,i'_b} q^{2,3}_{i''_1,j'_1}\cdots q^{2,3}_{i''_c,j'_c} q^{3,3}_{j''_1,k''_1}\cdots q^{3,3}_{j''_e,k''_e}.$$
In this notation, $b+c=r$, $c+e = s$, and $a+e = t$, and this lists $(i'_1,\dots,i'_b,i''_1,\dots,i''_c)$, $(j'_1,\dots,j'_c,j''_1,\dots,j''_e)$, and $(k'_1,\dots,k'_a,k''_1,\dots,k''_e)$ are permutations of the lists $(i_1,\dots,i_r)$, $(j_1,\dots,j_s)$, and $(k_1,\dots,k_t)$, respectively. 

Hence in the new variables $q^{a,b}_{i,j}, \tau^u_k$, the only possible monomials which can contain $\epsilon$ are of the form
\begin{equation} \tau^x_{k'_1}\cdots \tau^x_{k'_a} \tau^h_{i'_1}\cdots \tau^h_{i'_b} {\bf q}^{2,3}_{i''_1,j'_1}\cdots {\bf q}^{2,3}_{i''_c,j'_c} q^{3,3}_{j''_1,k''_1}\cdots q^{3,3}_{j''_e,k''_e}.\end{equation}
In this notation, ${\bf q}^{2,3}_{a,b}$ can denote either $q^{2,3}_{a,b}$ or $\tau^{y}_{a+b}$ (which contains $q^{2,3}_{a,b}$ by (4.16)). In order for $\epsilon$ to appear in $\omega$, any representative for the coset $\Phi^{-1}(\omega)\in F/I$ must contain at least one monomial 
\begin{equation}\mu=T^x_{k'_1}\cdots T^x_{k'_a} T^h_{i'_1}\cdots T^h_{i'_b} {\bf Q}^{2,3}_{i''_1,j'_1}\cdots {\bf Q}^{2,3}_{i''_c,j'_c} Q^{3,3}_{j''_1,k''_1}\cdots Q^{3,3}_{j''_e,k''_e}\end{equation} corresponding to (4.27). Here ${\bf Q}^{2,3}_{a,b}$ can denote either $Q^{2,3}_{a,b}$ or $T^{y}_{a+b}$, as above. In particular, the normal form $\hat{\omega}$ is a representative for $\Phi^{-1}(\omega)$, so it must contain at least one monomial $\mu$ of the form (4.28). Since every monomial appearing in $\hat{\omega}$ is divisible by $T^h_0$, it follows that $b>0$, and at least one of elements in the list $(i'_1,\dots,i'_b)$ is zero. We may assume without loss of generality that $i'_1 = 0$.

First, we claim that $e=0$. By Lemma 4.11, the factor $T^h_{0}Q^{3,3}_{j''_u k''_u}$ is the leading term of an element of the Grobner basis $B$, for each $u=1,\dots,e$. Since $\mu$ is in normal form and is divisible by $T^h_0$, it cannot be divisible by any of the factors $Q^{3,3}_{j''_u k''_u}$, so we have $e=0$. Hence $\mu$ has the form
$$T^x_{k_1}\cdots T^x_{k_t} T^h_0 T^h_{i'_2}\cdots T^h_{i'_b} {\bf Q}^{2,3}_{i''_1,j_1}\cdots {\bf Q}^{2,3}_{i''_s,j_s} .$$ Finally, we claim that for each $u=1,\dots,s$, we have ${\bf Q}^{2,3}_{i''_u,j_u} = T^y_{i''_u+j_u}$. Otherwise, $\mu$ would be divisible by $T^h_0 Q^{2,3}_{i''_u,j_u}$ for some $i''_u>0$, which is impossible since this is the leading term of an element of $B$, by  Lemma 4.11. Hence the monomial 
$$T^x_{k_1}\cdots T^x_{k_t} T^h_0 T^h_{i'_2}\cdots T^h_{i'_b}T^y_{i''_1+j_1}\cdots T^y_{i''_s+j_s} $$ appears in $\hat{\omega}$, as desired.

\end{document}